\documentclass[11pt,reqno]{amsart}

\newcommand{\mysection}[1]{\section{#1}
      \setcounter{equation}{0}}
 
\usepackage[mathscr]{euscript}

\newtheorem{theorem}{Theorem}[section]
\newtheorem{lemma}[theorem]{Lemma}

\newtheorem{corollary}[theorem]{Corollary}

\theoremstyle{definition}
\newtheorem{assumption}{Assumption}[section]
\newtheorem{definition}{Definition}[section]
\theoremstyle{remark}
\newtheorem{remark}{Remark}[section]

\newcommand\cbrk{\text{$]$\kern-.15em$]$}} 
\newcommand\opar{\text{\raise.2ex\hbox{${\scriptstyle | }$}\kern-.34em$($} }
\newcommand{\tr}{\text{\rm tr}\,}

 \makeatletter
 \def\dashint{%
 \operatorname%
 {\,\,\text{\bf--}\kern-.98em\DOTSI\intop\ilimits@\!\!}}
 \makeatother

\newcommand{\WO}{\overset{\scriptscriptstyle0}%
{W}\,\!}

\newcommand\gb{\mathfrak{b}}

\newcommand\bR{\mathbb{R}}

\newcommand\bL{\mathbb{L}}

\newcommand\bW{\mathbb{W}}

\newcommand\cB{\mathcal{B}}
\newcommand\cF{\mathcal{F}}

\newcommand\cP{\mathcal{P}}

\newcommand\cL{\mathcal{L}}

\newcommand\cW{\mathcal{W}}

\newcommand\frD{\mathfrak{D}}

\newcommand\sfB{{\textsc b}}

\begin{document}

\title[Kalman-Bucy filter and SPDEs]{
Kalman-Bucy filter and SPDEs
with growing lower-order coefficients in $W^{1}_{p}$
spaces without weights}

\author[N.V.  Krylov]{N.V. Krylov}%
\thanks{The work   was partially supported
  by NSF grant DMS-0653121}
\address{127 Vincent Hall, University of Minnesota,
Minneapolis,
       MN, 55455, USA}
\email{krylov@math.umn.edu}

\dedicatory{Dedicated to D.L. Burkholder}


\begin{abstract}
We consider divergence form uniformly parabolic SPDEs with
 VMO bounded leading coefficients, bounded  
coefficients in the stochastic part,
 and possibly growing lower-order coefficients in the deterministic part.
 We look for solutions
  which are summable to the $p$th power, $p\geq2$,
 with respect to the usual
Lebesgue measure along with their first-order   derivatives
with respect to the spatial variable.

Our methods allow us to include Zakai's equation
for the Kalman-Bucy filter into the general filtering theory. 
\end{abstract}

\maketitle

\mysection{Introduction}
We consider divergence form uniformly parabolic SPDEs with
bound\-ed VMO  leading coefficients, bounded  
coefficients in the stochastic part,
 and possibly growing lower-order coefficients in the deterministic part.
 We look for solutions
  which are summable to the $p$th power, $p\geq2$,
 with respect to the usual
Lebesgue measure along with their first-order  derivatives
with respect to the spatial variable. 
The present paper seems to be the first
one treating the unique solvability of
these equations  
without imposing any {\em special\/}
 conditions on the relations between
the coefficients or on their {\em derivatives\/}.

This article in its spirit
is similar to the author's recent articles \cite{KP},
\cite{Kr09_1},  \cite{Kr10_1}, and \cite{Kr10_2}
and we spare the reader the common part of the
comments about the literature,
which can be found in the above references. 
The main idea, we use, originated from \cite{KP} 
and \cite{Kr09_1} and relies on application of special
 cut-off functions whose support evolves in time
in a manner adapted to the drift terms.
The paper consists of two parts: Sections 
\ref{section 11.11.1} trough \ref{section 1.1.2}
are devoted to some general issues of the theory
of SPDEs with growing coefficients and in
Sections \ref{section 11.8.1} through
\ref{section 1.1.4} we apply the results of the
previous sections to show that the filtering equations
corresponding to the Kalman-Bucy filter fall into
the general theory.

In a sense
the methods of the 
first part of the present article arose as a combination of 
the methods from \cite{Kr10_1} and \cite{Kr10_2} which allow us
to combine the method used for PDE equations with irregular (VMO)
higher-order coefficients, growing lower-order
coefficients, and $p>1$ with the methods which work
in similar situation for SPDEs if $p=2$. 
Since we are interested in higher regularity of solutions
(see, for instance, Theorem \ref{theorem 1.12.1})
we use the power of summability $p\geq2$ and, in contrast
with \cite{Kr10_1}, this forces us to
require some regularity of the higher-order coefficients.
Roughly speaking we need the second-order
coefficients of the deterministic part of the equation
belong to VMO in $x$ and the first-order coefficients
of the stochastic part to be uniformly continuous in $x$.
In particular, the results of the present article
do not generalize those of \cite{Kr10_1}.

On the other hand,
if we drop all stochastic terms, then we obtain
the results of \cite{Kr10_2} for $p\geq2$, which by duality,
available for deterministic equations, allows one to
extend the result to full range $p>1$.
Concerning the {\em deterministic\/}
 equations with growing coefficients
in spaces with or without weights
it is worth mentioning that

(i) Equations in spaces with weights are treated, for instance, in
\cite{CV1}, \cite{ChG}, \cite{FL}, \cite{Lu}, and
 \cite{MP1} for time independent coefficients, part of the result of which
are extended in \cite{GL} to
time-dependent Ornstein-Uhlenbeck operators;

(ii) Equations in spaces without weights are treated, 
for instance, in \cite{LV}, \cite{MP}, \cite{PR},
and \cite{CF}.

Some conclusions in the above cited papers are quite similar 
to ours
but the corresponding assumptions are not as general
in what concerns the regularity of the coefficients.
However, these papers contain a lot of additional
important information, which is probably
impossible to obtain by using our methods.

The second part of the article is devoted to the Kalman-Bucy
filter. One can say that one of the   sources
of interest in
SPDEs with growing coefficients is Zakai's equation
for filtering density in the case of
 partially observable diffusion processes. This equation
has divergence form which makes it possible to use
the results of the first part of the article.
In a very particular case of Gaussian processes
the filtering density is given by the
 Kalman-Bucy filter. 
 Generally, part of the 
coefficients of filtering  equations in 
case of Gaussian processes grow.  When the 
coefficients of an SPDE 
grow, it is quite natural to consider the equations
in function 
 spaces with weights which would restrict the set of solutions
in such a way that all terms in the equation will be
from the same space as the free terms. There are very many articles
which use this idea in $\cL_{2}$- and $\cL_{p}$-settings
(see, for instance, \cite{CV}, \cite{GK}, \cite{Gy93},
\cite{Gy97} and the references therein). Unfortunately, 
the application
of the spaces with weights do not allow one
 to treat filtering equations
corresponding to the Kalman-Bucy filter even without the
 so-called cross terms when the operators $\Lambda^{k}_{t}$
in \eqref{10.20.2} are of zeroth order. The main obstacle here is
that the zeroth order coefficient of $\Lambda^{k}_{t}$
is a linear function of $x$. In the general theory, which
we develop  in this article, we do not allow it to grow
either and we use an auxiliary function to ``kill'' this coefficient.
The construction of this auxiliary function exploits
a specific structure of the equation and allows us
to transform the general filtering
equation \eqref{10.20.2} to its ``reduced"
 form \eqref{9.29.01}, which does not
contain the zeroth order term in the stochastic part.
After that one can use a simple change of the unknown function
shifting the $x$ variables in such a way that 
the stochastic part of \eqref{9.29.01} will disappear altogether
and the equation will become a parabolic equation
with time inhomogeneous and random Ornstein-Uhlenbeck
operator. The fact that the operator is time inhomogeneous
makes it impossible to apply any results based, for instance, on
the semigroup approach and  even specifically aimed at the
Ornstein-Uhlenbeck operator, which one can find in the 
above mentioned recent
articles such as \cite{ChG}, \cite{FL}, \cite{MP1},
or other results on elliptic operators with unbounded
coefficients such as in
\cite{PR}. The results of \cite{CV} 
are not applicable either because in \cite{CV}
the zeroth-order coefficient is assumed to grow quadratically
if the firs-order coefficients grow linearly.
However, the results of \cite{GK} on general SPDEs
with growing coefficients are applicable 
to the reduced form of the SPDE for the Kalman-Bucy filter
and they provide existence and uniqueness theorems
in Sobolev spaces with $p=2$ and weights depending on $t,x$
{\em and\/} $\omega$. By the way, a drawback of using
weights depending on $t$  is that
one cannot extract from the results for general SPDEs
any result for deterministic elliptic equations. 
 
 If one concentrates on $p=2$, then one can use the results from 
\cite{GL} where the Ornstein-Uhlenbeck time inhomogeneous
operators are investigated in Sobolev spaces with
Gaussian time dependent weight. Again this would
allow one to investigate \eqref{9.29.01} 
in Sobolev spaces with $p=2$ and weights depending on $t,x$
and $\omega$.  We deal with any $p\geq2$ and
do not use weights.

The article is organized as follows. In 
Section \ref{section 11.11.1} we introduce basic notation,
function spaces, and equations.
Section \ref{section 2.15.1} contains our main results
concerning SPDEs. Section \ref{section 1.1.1} contains
the proof of Theorem \ref{theorem 3.11.1} concerning an
apriori estimate and Theorem \ref{theorem 1.12.1}
about regularity properties of solutions. In Section \ref{section 6.9.5}
we prove the existence Theorem  \ref{theorem 3.16.1}.

In Section \ref{section 1.1.2} we prove a version of It\^o's
formula which allows us to use the results of the previous
sections to derive the filtering equation without
using anything from the filtering theory itself.
We do it by following \cite{KZ} and \cite{Kr_10}.
In Section \ref{section 11.8.1} we state our main result
about the equation corresponding to Kalman-Bucy filter.
We consider the so-called conditionally Gaussian process
in the spirit of \cite{LS}. However, in contrast 
with \cite{LS}, our coefficients
depend only on the current state of the two-component
process under consideration and are
not allowed to depend on the whole past 
of the observable component.
In Section \ref{section 1.1.3} we consider
the ``reduced'' form \eqref{9.29.01}
of the main filtering equation
\eqref{10.20.2}. The results of the previous sections
turn out to be applicable to \eqref{9.29.01}.
In  the final Section \ref{section 1.1.4}
we finish proving
Theorems \ref{theorem 11.8.1}
and  \ref{theorem 11.11.1}, part of assertions of the former
being proved in Section~\ref{section 1.1.3}.

 \mysection{General setting}
                                      \label{section 11.11.1}

Let $(\Omega,\cF,P)$ be a complete probability space
with an increasing filtration $\{\cF_{t},t\geq0\}$
of complete with respect to $(\cF,P)$ $\sigma$-fields
$\cF_{t}\subset\cF$. Denote by $\cP=\cP(\{\cF_{t}\})$ 
the predictable
$\sigma$-field in $\Omega\times(0,\infty)$
associated with $\{\cF_{t}\}$. Let
 $w^{k}_{t}$, $k=1,2,...$, be independent one-dimensional
Wiener processes with respect to $\{\cF_{t}\}$.  
Let $\tau$ be a stopping time.

We consider the   second-order
operator $L_{t}$
\begin{equation}                                    \label{lu}
 L_{t} u_{t}(x) =D_{i}\big( a^{ij}_{t}( x)D_{j}u_{t}(x)
+\gb^{i}_{t}(x)u_{t}(x)\big) + b^{i} _{t}(x)
D_{i} u_{t}(x) -c _{t}(x) u_{t}(x),
\end{equation}
and the first-order operators
$$
\Lambda^{k}_{t} u_{t}(x)=\sigma^{ik}_{t}(x)D_{i}u_{t}(x)
+\nu^{k}_{t}(x)u_{t}(x)
$$  
   acting on functions $u_{t}(x)$ defined
on $\Omega\times\bR^{d+1}_{+}$, where $\bR^{d+1}_{+}=
 [0, \infty) \times \bR^d$, and given for $k=1,2,...$ 
(the summation convention is enforced throughout the article),
where  
$$
D_{i}=\frac{\partial}{\partial x^{i}}.
$$
We set  $\bR_{+}=[0,\infty)$.

Our main concern in the first part of the paper
is proving the unique solvability
of the equation
\begin{equation}
                                                \label{2.6.4}
du_{t}=(L_{t}u_{t}-\lambda u_{t}+D_{i}f^{i}_{t}+f^{0}_{t})\,dt
+(\Lambda^{k}_{t}u_{t}+g^{k}_{t})\,dw^{k}_{t},
\quad t\leq\tau,
\end{equation}
with an appropriate initial condition at $t=0$, where  
$\lambda\geq0$ is a constant. 
The precise 
assumptions on the coefficients, free terms, and initial data
will be given later. First we introduce appropriate function
spaces.

Fix a number 
$$
p\geq2, 
$$
and denote $\cL_{p}=\cL_{p}(\bR^{d})$.
We use the same notation $\cL_{p}$ for vector- and matrix-valued
or else
$\ell_{2}$-valued functions such as
$g_{t}=(g^{k}_{t})$ in \eqref{2.6.4}. For instance,
if $u(x)=(u^{1}(x),u^{2}(x),...)$ is 
an $\ell_{2}$-valued measurable function on $\bR^{d}$, then
$$
\|u\|^{p}_{\cL_{p}}=\int_{\bR^{d}}|u(x)|_{\ell_{2}}^{p}
\,dx
=\int_{\bR^{d}}\big(
\sum_{k=1}^{\infty}|u^{k}(x)|^{2}\big)^{p/2}
\,dx.
$$
 
As usual,  
$$
W^{1}_{p}=\{u\in \cL_{p}: Du\in \cL_{p}\},
\quad
 \|u\|_{W^{1}_{p}}=
\|u\|_{\cL_{p}}+\|Du\|_{\cL_{p}},
$$
where by $Du$   we mean the gradient  with respect
to $x$ of a function $u$ on $\bR^{d}$.

Recall that $\tau$ is a stopping time and introduce
$$
\bL _{p}(\tau):=\bL _{p}(\{\cF_{t}\},\tau):=\cL_{p}(\opar 0,\tau\cbrk,\cP,
\cL_{p}),
$$
$$
\bW^{1}_{p}(\tau):=\bW^{1}_{p}(\{\cF_{t}\},\tau):=
\cL_{p}(\opar 0,\tau\cbrk,\cP,
W^{1}_{p}),
$$
$$
\bL_{p}=\bL_{p}(\infty),
\quad\bW^{1}_{p}=\bW^{1}_{p}(\infty).
$$
Remember that the elements of
$\bL_{p}( \tau)$ need only 
belong to $\cL_{p}$ on a predictable subset of 
$\opar 0,\tau\cbrk$ of full measure. For the sake of convenience
we will always assume that they are defined everywhere
on $\opar 0,\tau\cbrk$ at least as generalized functions. 
Similar situation occurs in the case of $\bW^{1}_{p}( \tau)$.

The following definition is most appropriate
for investigating our equations
if the coefficients of $L_{t}$ and $\Lambda^{k}_{t}$
are bounded.
\begin{definition}
                                    \label{definition 3.16.1}
 Introduce $\cW^{1}_{p}(\tau)$,
as the space of functions $u_{t}
=u_{t}(\omega,\cdot)$ on $\{(\omega,t):
0\leq t\leq\tau(\omega),t<\infty\}$ with values
in the space of generalized functions on $\bR^{d}$
and having the following properties:

(i) We have $u_{0}\in \cL_{p}(\Omega,\cF_{0},\cL_{p})$;

(ii)  We have $u
\in \bW^{1}_{p}(\tau )$;

(iii) There exist   $f^{i}\in \bL_{p}(\tau)$,
$i=0,...,d$, and $g=(g^{1},g^{2},...)\in \bL_{p}(\tau)$
such that
 for any $\varphi\in C^{\infty}_{0}=C^{\infty}_{0}(\bR^{d})$ 
with probability 1
for all   $t\in[0,\infty)$
we have
$$
(u_{t\wedge\tau},\varphi)=(u_{0},\varphi)
+\sum_{k=1}^{\infty}\int_{0}^{t}I_{s\leq\tau}
(g^{k}_{s},\varphi)\,dw^{k}_{s}
$$
\begin{equation}
                                                 \label{1.2.1}
+
\int_{0}^{t}I_{s\leq\tau}\big((f^{0}_{s},\varphi)-(f^{i}_{s},D_{i}\varphi)
 \big)\,ds.
\end{equation}
In particular, for any $\phi\in C^{\infty}_{0}$, the process
$(u_{t\wedge\tau},\phi)$ is $\cF_{t}$-adapted and (a.s.) continuous.
In case that property (iii) holds, we write
$$
du_{t}=(D_{i}f^{i}_{t}+f^{0}_{t})\,dt+g^{k}_{t}\,dw^{k}_{t},
\quad t\leq\tau.
$$
Finally, set $\cW^{1}_{p}=\cW^{1}_{p}(\infty)$.
\end{definition}

\begin{remark}
                                              \label{remark 1.12.3}
The reader understands that if $u$ is a generalized
function on $\bR^{d}$, then $(u,\phi)$ represents
the result of the action of $u$ on the test function
$\phi\in C^{\infty}_{0}$. When $u$ is a locally integrable
function, $(u,\phi)$ is the integral of the product $u\phi$.
According to these notation
$$
(f^{0}_{s},\varphi)-(f^{i}_{s},D_{i}\varphi)=
(\bar{f}_{s},\phi),
$$
where the function $\bar{f}_{s}$ with values in the space
of generalized functions is defined by $\bar{f}_{s}
=D_{i}f^{i}_{s}+f^{0}_{s}$. In the framework of Definition
\ref{definition 3.16.1} we have $\bar{f}
\in\cL_{p}(\opar0,\tau\cbrk,\cP,H^{-1}_{p})$, where
$H^{-1}_{p}=(1-\Delta)^{1/2}\cL_{p}$. One also
knows that any $\bar{f}
\in\cL_{p}(\opar0,\tau\cbrk,\cP,H^{-1}_{p})$ is written
as $\bar{f}_{s}
=D_{i}f^{i}_{s}+f^{0}_{s}$ with some $f^{j}\in\bL_{p}(\tau)$.
\end{remark}

Also introduce
the spaces of initial data in the same way as in \cite{Kr09}.
\begin{definition}
                                        \label{definition 8.10.1}
Let $u_{0}$ be an $\cF_{0}$-measurable function on $\Omega$
with values in the space of generalized functions
on $\bR^{d}$. We write
$u_{0}\in \tr \cW^{1}_{p}=\tr\cW^{1}_{p}(\cF_{0})$ 
if there exists a function 
$v \in   \cW^{1}_{p} 
 $ such that $d v_{t}=(\Delta v_{t}-v_{t})\,dt$, $t\in\bR_{+}$,
and $v_{0}=u_{0}$. In such a case we set
$$
\|u_{0}\|^{p}_{\tr \cW^{1}_{p}}=E\|v\|^{p}_{\bW^{1}_{p} }.
$$
\end{definition}

One knows that $\tr \cW^{1}_{p}$ is a Banach space,
$v$ in the above definition is unique and $\cF_{0}$-measurable.

We  give the definition of
solution of  \eqref{2.6.4} adopted throughout the article
and which in case the coefficients of 
$L_{t}$ and $\Lambda^{k}_{t}$
are bounded coincides with the one obtained by applying 
 Definition \ref{definition 3.16.1}.
 
\begin{definition}
                                  \label{definition 3.20.01}
Let $f^{j}\in\bL_{p}(\tau)$, $j=0,...,d$,
$g=(g^{1},g^{2},...) \in\bL_{p}(\tau)$.
By a solution of
\eqref{2.6.4} (relative to $\{\cF_{t}\}$)  with initial condition
$u_{0}\in \tr \cW^{1}_{p}$
we mean a function $u\in\bW^{1}_{p}(\tau)$ (not
$\cW^{1}_{p}(\tau)$) such that

(i)  For any $\phi\in C^{\infty}_{0} $ 
 the integrals in
$$
(u_{t\wedge\tau},\phi)=(u_{0},\phi)
+\sum_{k=1}^{\infty}\int_{0}^{t}I_{s\leq\tau}
(\sigma^{ik}_{s}D_{i}u_{s}+\nu^{k}_{s}u_{s}
+g^{k}_{s},\phi)\,dw^{k}_{s}
$$
\begin{equation}
                                                \label{3.16.7}
+
\int_{0}^{t}I_{s\leq\tau}\big[(b^{i}_{s}D_{i}u_{s}
-(c_{s}+\lambda)u_{s}+f^{0}_{s},\phi)
-(a^{ij}_{s}D_{j}u_{s}+\gb^{i}_{s}u_{s}+
f^{i}_{s},D_{i}\phi)
 \big]\,ds
\end{equation}
are 
well defined and are finite for all finite
 $t\in\bR_{+}$ and the series converges uniformly
on finite subinterval of $\bR_{+}$ in probability;

(ii) For any $\phi\in C^{\infty}_{0} $ with probability one
equation  \eqref{3.16.7}  
 holds for all  $t\in\bR_{+}$.

\end{definition}

Observe that for any solution of  \eqref{2.6.4}
in the sense of the above definition and any 
$\phi\in C^{\infty}_{0}$ the process $(u_{t\wedge\tau},\phi)$
is continuous (a.s.) and $\cF_{t}$-adapted.

Also notice that, if the coefficients of $L$ and $\Lambda^{k}$
are bounded, then any $u\in\cW^{1}_{p}(\tau)$ is a solution
of \eqref{2.6.4} with appropriate free terms since if
\eqref{1.2.1} holds, then \eqref{2.6.4} holds 
(always in the sense of Definition \ref{definition 3.20.01})
as well with
$$
f^{i}_{t}-a^{ij}_{t}D_{j}u_{t}-\gb^{i}u_{t},
\quad i=1,...,d,\quad
f^{0}_{t}+(c_{t}+\lambda)u_{t}-b^{i}_{t}D_{i}u_{t},
$$
$$
g^{k}_{t}-\sigma^{ik}D_{i}u_{t}-\nu^{k}_{t}u_{t}
$$
in place of $f^{i}_{t}$, $i=1,...,d$, $f^{0}_{t}$, and $g^{k}_{t}$,
respectively.

\mysection{Main results for SPDEs}
                                   \label{section 2.15.1}
For $\rho>0$ denote $B_{\rho}(x)=\{y\in\bR^{d}:|x-y|<\rho\}$,
$B_{\rho}=B_{\rho}(0)$. 
\begin{assumption} 
                                       \label{assumption 2.7.2}
(i) The functions $a^{ij}_{t}(x)$, $\gb^{i}_{t}(x)$,
$b^{i}_{t}(x)$, $c_{t}(x)$,
$\sigma^{ik}_{t}(x)$, $\nu^{k}_{t}(x)$ are real valued, measurable
with respect to $\cF\otimes\cB(\bR^{d+1}_{+})$,
 $\cF_{t}$-adapted for any $x$, and  $c\geq 0$.

(ii) There exists a constant  $ \delta>0$ such that
for all values of arguments and $\xi\in\bR^{d}$
$$
(a^{ij}   
-  \alpha^{ij} ) \xi^{i}
\xi^{j}\geq\delta|\xi|^{2},\quad
|a^{ij}|\leq \delta^{-1} ,
\quad  |\nu|_{\ell_{2}}\leq \delta^{-1},
$$
where  $\alpha^{ij} =(1/2)(\sigma^{i\cdot},\sigma^{j\cdot})
_{\ell_{2}}$. Also, the constant $\lambda\geq0$.

(iii) For any  $x\in
\bR^{d}$ (and $\omega$) the function
\begin{equation}
                                                \label{9.3.01}
\int_{B_{1}}(|\gb_{t}(x+y)|+|b_{t}(x+y)|+|c_{t}(x+y)|)\,dy
\end{equation}
is locally   integrable to the $p'$th power on $\bR_{+}=[0,\infty)$,
where $p'=p/(p-1)$.
\end{assumption}

Notice that the matrix $a=(a^{ij})$ need not be symmetric.
Also notice that in Assumption \ref{assumption 2.7.2} (iii)
the ball $B_{1}$ can be replaced with any other ball
without changing the set of admissible coefficients
$\gb,b,c$.

Recall that as is well known
  if $u\in\bW^{1}_{p}(\tau)$, then owing to the boundedness of $\nu$
and $\sigma$ and the fact that $Du,u,g\in\bL_{p}(\tau)$, $p\geq2$,
the first series on the right in \eqref{3.16.7}
converges uniformly in probability and the series
is a continuous local martingale. Furthermore, if we denote it by
$m_{t}$, then for  any $T\in\bR_{+}$
$$
E\sup_{t\leq T}|m_{t}|^{p}\leq N
E\big(\sum_{k=1}^{\infty}\int_{0}^{\tau\wedge T} 
(\sigma^{ik}_{s}D_{i}u_{s}+\nu^{k}_{s}u_{s}
+g^{k}_{s},\phi)^{2}\,ds\big)^{p/2}
$$
$$
\leq N\|\phi\| _{\cL_{1}}^{p/2}
E\big( \int_{0}^{\tau\wedge T}\sum_{k=1}^{\infty}
(|\sigma_{s}^{ik}|^{2}|D_{i}u_{s}|^{2}+|\nu^{k}_{s}|^{2}
|u_{s}|^{2}
+|g^{k} _{s}|^{2},|\phi|) \,ds\big)^{p/2}
$$
\begin{equation}
                                                    \label{9.3.3}
\leq N(\|u\|^{p}_{\bW^{1}_{p}(\tau)}+
\|g\|^{p}_{\bL_{p}(\tau)}),
\end{equation}
where the constants $N$ depend only on $\phi$, $d$,
$p$, $\delta$, and $T$.

 \begin{assumption} 
                                       \label{assumption 2.4.1}
There exists a function $\kappa(r)$,
$r\in\bR_{+}$, such that $\kappa(0+)=0$ and for
 any $\omega\in\Omega$, 
 $t\geq0$, $x,y\in\bR^{d}$,
 and $i =1,...,d$  we have 
$$
|\sigma_{t}^{i\cdot}(x)-\sigma^{i\cdot}_{t}(y)|_{\ell_{2}}
\leq\kappa(|x-y|).
$$

\end{assumption}

The following assumptions contain parameters   $ 
\gamma_{a},\gamma_{b}  \in(0,1]$,  
whose values  will be specified later.
They also contain  constants $K\geq0$, $\rho_{0}, \rho_{1}
\in(0,1]$
 which are fixed.

 \begin{assumption} 
                                       \label{assumption 2.4.01}
For any $\omega\in\Omega$, 
$\rho\in(0,\rho_{0}]$,  $t\geq0$,
 and $i,j=1,...,d$  we have 
\begin{equation}
                                             \label{2.9.5}
\rho^{-2d-2}\int_{t}^{t+\rho^{2}}\bigg(
 \sup_{x\in \bR^{d}}\int_{B_{\rho}(x)}\int_{B_{\rho}(x)}
|a^{ij}_{s}(  y)-a^{ij}_{s}( z)|\,dydz\bigg)
 \,ds\leq\gamma_{a}.
\end{equation}

\end{assumption}

Obviously, the left-hand side of \eqref{2.9.5}
is less than
$$
N(d)\sup_{t\geq0}\sup_{|x-y|\leq2\rho}|a^{ij}_{t}( x )
-a^{ij}_{t}(y)|,
$$
which implies that Assumption \ref{assumption 2.4.01}
 is satisfied with any $\gamma_{a} >0$ if, for instance,
$a$ is   uniformly continuous in $x$ uniformly in $\omega$ and $t$.
Recall that if $a$ is independent of $t$ and
 for any $\gamma_{a} >0$ there is a $\rho_{0}>0$
such that Assumption \ref{assumption 2.4.01}
 is satisfied, then one says that
$a$ is in VMO.

We take and fix a number $q=q(d,p) $ such  that
\begin{equation}
                                                  \label{8.11.3}
q\geq\max(d,p)\quad\hbox{\rm if}\quad p\ne d,
\quad q>d\quad\hbox{\rm if}\quad p= d.
\end{equation}

\begin{assumption} 
                                      \label{assumption 3.11.1} 
For any $\omega\in\Omega$,
 $ \gb:=(\gb^{1} ,...,\gb^{d} )$,
  $ b:=(b^{1} ,...,b^{d} )$, and   $(t,x)\in\bR^{d+1}$
 we have  
$$
 \int_{B_{\rho_{1}}(x)}\int_{B_{\rho_{1}}(x)}|\gb_{t}( y)
-\gb_{t}(   z)|^{q }\,dydz  +
 \int_{B_{\rho_{1}}(x)}\int_{B_{\rho_{1}}(x)}|b_{t}( y)
-b_{t}(  z)|^{q}\,dydz   
$$
$$
+ \int_{B_{\rho_{1}}(x)}\int_{B_{\rho_{1}}(x)}|c_{t}( y)
-c_{t}( z)|^{q}\,dydz  \leq KI_{q>d}+\rho_{1}^{d}\gamma_{b}.
$$

\end{assumption}

Obviously,  Assumption \ref{assumption 3.11.1}  
is satisfied   if  
 $b$, $\gb$, and $c$  are independent of $x$.
They also are satisfied  with any $q>d$, $\gamma_{b}=0$,
and $\rho_{1}=1$ on the 
account of choosing
$K$ appropriately if, say,
$$
|\gb_{t}(  x )
-\gb_{t}(y)|+|b_{t}(  x )
-b_{t}( y)|+|c_{t}(  x )
-c_{t}(  y)|\leq K_{1}
$$ 
whenever $|x-y|\leq1$, where $K_{1}$ is a constant.
In particular,  Assumption \ref{assumption 3.11.1}  
is satisfied   if 
$\gb$, $b$, and $c$ are globally Lipschitz continuous:
$$
|\gb_{t}(  x )
-\gb_{t}(y)|+|b_{t}(  x )
-b_{t}( y)|
$$
\begin{equation}
                                                   \label{1.15.1}
+|c_{t}(  x )
-c_{t}(  y)|\leq K_{1}|x-y|,\quad\forall x,y\in\bR^{d},t\geq0.
\end{equation} 
 We see   that Assumption \ref{assumption 3.11.1}
 allows $b$, $\gb$, and $c$ growing linearly in $x$.
Here is our result on apriori estimates of solutions of
\eqref{2.6.4}.
 
\begin{theorem}
                                       \label{theorem 3.11.1}
There exist
$$
\gamma_{a} =\gamma_{a}(d,\delta,p),\quad
\gamma_{b} =\gamma_{b}(d,\delta,p,\kappa,\rho_{0} )\in(0,1],
$$
$$
 N=N(d,\delta,p,\kappa,\rho_{0}   ),
\quad \lambda_{0}=\lambda_{0}(d,\delta, p,\kappa,
\rho_{0},\rho_{1},K)\geq1
$$  
such that, if the above assumptions are satisfied
and $\lambda\geq
\lambda_{0}$ and
  $u  $ is a solution of \eqref{2.6.4}
with   initial data $u_{0}\in \tr \cW^{1}_{p}$ and some
$f^{j},g \in\bL_{p}(\tau)$, then
$$
\lambda\|u\|^{2}_{\bL _{p}(\tau)}+\|Du\|^{2}_{\bL _{p}(\tau)}
\leq N\big(\sum_{i=1}^{d}\|f^{i}\|^{2}_{\bL _{p}(\tau)}
+\|g\|^{2}_{\bL _{p}(\tau)} \big)
$$
\begin{equation}
                                       \label{3.11.2}
+N\lambda^{-1}\|f^{0}\|^{2}_{\bL _{p}(\tau)}
+N \|u_{0}\|^{2}_{\tr \cW^{1}_{p}} .
\end{equation}
\end{theorem}

\begin{remark}
                                      \label{remark 9.27.1}
There is an unusual property of $u_{t}$,
which is nontrivial even if $ f^{j}_{t}
=g^{k}_{t}
\equiv0$. 

Namely, assume that   $g\equiv0$. 
Take a predictable $\ell_{2}$-valued process
$\xi_{t}$ such that $(\nu_{t},\xi_{t})_{\ell_{2}}\geq0$ 
and $(\nu_{t},\xi_{t})_{\ell_{2}}$ and $(\sigma^{i\cdot}_{t},
\xi_{t})$ are independent of $x$ (which happens, for
instance, if $\nu=0$ and $\sigma$ is independent of $x$)
and
$$
\int_{0}^{\tau}|\xi_{t}|_{\ell_{2}}^{2}\,dt<\infty
$$
(a.s.) and assume that 
$E\rho_{\tau}(\xi)=1$, where
$$
\rho_{t}(\xi)=\rho_{t}(\xi,dw):=
\exp\big(-\int_{0}^{t}\xi^{k}_{s}\,dw^{k}_{s}
-\tfrac{1}{2}\int_{0}^{t}|\xi_{s}|_{\ell_{2}}^{2}\,ds\big) .
$$
Then the assertion of Theorem \ref{theorem 3.11.1}
holds with the {\em same\/} $\gamma_{a}$,
$\gamma_{b}$, $\lambda_{0}$, and $N$
if we understand $\|v\|_{\bL_{p}(\tau)}^{p}$ for all $v$'s as  
$$
E\rho_{\tau}\int_{0}^{\tau}\|v_{t}\|^{p}_{\cL_{p}}\,dt.
$$

Indeed, one can change
the probability measure by using Girsanov's theorem.
This will add a new drift term in the deterministic part
of  \eqref{2.6.4}  and this  additional drift depends only on
$(\omega, t)$. This will also add the
term $-(\nu_{t},\xi_{t})_{\ell_{2}}u_{t}\,dt$, 
where $(\nu_{t},\xi_{t})_{\ell_{2}}$ is
 nonnegative and also independent
of $x$. Then the result follows
immediately from 
Theorem~\ref{theorem 3.16.1}.  

\end{remark}
Theorem \ref{theorem 3.11.1} admits the following version
if $\tau$ is bounded.
\begin{theorem}
                                       \label{theorem 1.14.1}
Let $T\in(0,\infty)$ be a constant and 
suppose that $\tau\leq T$.
Assume that the above assumptions are satisfied
with $\gamma_{a}$ and $\gamma_{b}$ from
 Theorem~\ref{theorem 3.11.1}. Let $\lambda=0$ and let
  $u  $ be a solution of \eqref{2.6.4}
with   initial data $u_{0}\in \tr \cW^{1}_{p}$ and some
$f^{j},g \in\bL_{p}(\tau)$. Then
\begin{equation}
                                       \label{1.14.2}
\|u\|^{2}_{\bW^{1} _{p}(\tau)}
\leq N\big(\sum_{i=0}^{d}\|f^{i}\|^{2}_{\bL _{p}(\tau)}
+\|g\|^{2}_{\bL _{p}(\tau)} 
+\|u_{0}\|^{2}_{\tr \cW^{1}_{p}}\big) ,
\end{equation}
where $N=N(d,\delta, p,\kappa,
\rho_{0},\rho_{1},K,T)$.
\end{theorem}

This result is a trivial consequence of Theorem
\ref{theorem 3.11.1} since, for any constant $\mu$, the function
$v_{t}:=u_{t}e^{-\mu t}$ satisfies \eqref{2.6.4} with
$\lambda+\mu$, $f^{j}_{t}e^{-\mu t}$, and $g^{k}_{t}e^{-\mu t}$
in place of $\lambda$, $f^{j}_{t}$, and $g^{k}_{t}$,
respectively. If $\mu$ is large enough and $\tau\leq T$, estimate
\eqref{3.11.2} for $v$ implies \eqref{1.14.2} indeed.

\begin{remark}
Theorems \ref{theorem 3.11.1} and \ref{theorem 1.14.1}
 provide uniqueness of solutions
of \eqref{2.6.4}. The apriori estimates \eqref{3.11.2}
and \eqref{1.14.2}
can also be used to investigate continuous dependence
of solutions on the coefficients and other data.

\end{remark}

To prove the existence we need stronger assumptions
because, generally, Assumption \ref{assumption 3.11.1}
does not guarantee that
$$
D_{i}(\gb^{i}_{t}u_{t})+b^{i}_{t}D_{i}u_{t}-c_{t}u_{t}
$$
can be written even locally as 
$D_{i}\hat{f}^{i}_{t}+\hat{f}^{0}_{t}$
with $\hat{f}^{j}\in\bL_{p}(\tau)$ if we only know that
$u\in\bW^{1}_{p}(\tau)$ even if $\gb$, $b$, and $c$ are independent
of $x$. 
We can only prove our   Lemma \ref{lemma 3.23.2}
 if we have a certain control on this expression.

\begin{assumption}
                                       \label{assumption 3.16.1}
For any  $x\in
\bR^{d}$ (and $\omega$) the function \eqref{9.3.01}  
is locally   integrable to the power
$ p /(p-2)$ (locally bounded if $p=2$) on $\bR_{+}=[0,\infty)$.
 \end{assumption}
\begin{remark}
                                            \label{remark 1.15.1}
Assumptions \ref{assumption 3.11.1} and \ref{assumption 3.16.1} 
are both satisfied if the global Lipschitz condition
\eqref{1.15.1} holds and $b_{t}(0)$, $\gb_{t}(0)$, and
$c_{t}(0)$ are bounded for each~$\omega$.
\end{remark}

\begin{theorem}
                                        \label{theorem 3.16.1}
Let the above assumptions be satisfied with
$\gamma_{a} $ and $\gamma_{b}$ taken from Theorem \ref{theorem 3.11.1}.
Take   $\lambda\geq\lambda_{0}$, where $\lambda_{0}$
is defined in Theorem \ref{theorem 3.11.1}, and take 
$u_{0}\in\tr \cW^{1}_{p}$. Then
there exists
a unique solution of \eqref{2.6.4}
  with initial condition $u_{0}$.

\end{theorem}

\begin{remark}
                                       \label{remark 9.29.2} 

If the stopping time $\tau$ is bounded, then in the above theorem
one can take $\lambda_{0}=0$. This is shown by the same argument
as after Theorem~\ref{theorem 1.14.1}.

\end{remark}

In general the continuity properties in $t$ of the solution
from Theorem \ref{theorem 3.16.1} are unknown.
For instance,
we do not know if $\|u_{t\wedge\tau}\phi\|_{\cL_{p}}$ is continuous
(a.s) for any $\phi\in C^{\infty}_{0}$. However, under
stronger assumptions we can say more about regularity
of $u$. In the following theorem by $H^{\gamma}_{p}$
we mean $(1-\Delta)^{-\gamma/2}\cL_{p}$.
\begin{theorem}
                                         \label{theorem 1.12.1}
Under the above assumptions suppose that
for each $x\in\bR^{d}$ the function \eqref{9.3.01} is
{\em bounded\/} on $\opar0,\tau\cbrk$. Then
the (unique) solution $u$ possesses the following properties:

(i) For any $\phi\in C^{\infty}_{0}$ we have
$\phi u\in\cW^{1}_{p}(\tau)$;

(ii) For any $\phi\in C^{\infty}_{0}$ the process
 $u_{t\wedge\tau}\phi$ is continuous on $\bR_{+}$
as an $\cL_{p}$-valued process (a.s.);

(iii) If $p>2$ and $\tau$ is bounded and we have two numbers
$\alpha$ and $\beta$ such that
$$
\frac{2}{p}<\alpha<\beta\leq1,
$$
 then
for any $\phi\in C^{\infty}_{0}$ (a.s.)
$$
u\phi\in C^{\alpha/2-1/p}([0,\tau], H^{1-\beta}_{p}).
$$
In particular, if $p>d+2$, then 

(a) for any $\varepsilon
\in(0,\varepsilon_{0}]$, with
$$
\varepsilon_{0}=1-\frac{d+2}{ p},$$
 (a.s.) for any
$t\in[0,\tau]$ we have $u_{t}\phi\in C^{\varepsilon_{0}-
\varepsilon}(\bR^{d})$ and the norm of $u_{t}\phi$
in this space is bounded as a function of $t$;

(b) for any $\varepsilon$ as in (a) (a.s.)
for  any $x\in\bR^{d}$
we have $u_{\cdot}(x)\phi(x)\in C^{(\varepsilon_{0}-
\varepsilon)/2}([0,\tau])$ and the norm of $ u_{\cdot}(x)\phi(x)$
in this space is bounded as a function of $x$.

\end{theorem}

Observe that assertions (ii) and (iii) 
of Theorem \ref{theorem 1.12.1} follow from
assertion (i) proved
in Remark \ref{remark 1.12.1}.
 In case of assertion (ii) this is shown
in \cite{Kr09_3}. The main part of assertion (iii) follows from
assertion (i) and Corollary 4.12 \cite{Kr01}. By applying
Sobolev's embedding theorems
 assertion (iii) (a) is obtained after taking $\alpha$
and $\beta$ close to $2/p$ and (iii) (b) after taking
$\alpha$
and $\beta$ close to $1-d/p$.

\begin{remark}
Let $p_{1},p_{2}\in[2,\infty)$,
let $\tau$ be bounded (cf. Remark \ref{remark 9.29.2}),
and let the assumptions of Theorem  
\ref{theorem 3.16.1} be satisfied for any $p\in[p_{1},p_{2}]$
 with   $\gamma_{a}$ and $\gamma_{b}$
 which are suitable for all 
$p\in[p_{1},p_{2}]$. Then it turns out  
that  the solution from Theorem  
\ref{theorem 3.16.1}  corresponding to
$p=p_{1}$ coincides with the one obtained for $p=p_{2}$.

This fact is obtained in the same way 
as the proof of Theorem 3.4 of \cite{Kr10_2}
is obtained from the proof of Theorem  3.3 \cite{Kr10_2}.
\end{remark}

Our last main result on general SPDEs bears on 
the measurability of $u_{t}$ with respect to 
$\sigma$-fields which are smaller than $\cF_{t}$.
It will be used in Section \ref{section 1.1.3}
 and this is the reason
why we use the somewhat strange notation $\tilde{y}_{t}$ and
$\tilde{\sfB}^{k}_{t}$
below.
We suppose that all
the above assumptions are satisfied with
$\gamma_{a}$ and $\gamma_{b}$
 taken from Theorem \ref{theorem 3.11.1}
and let $\tilde{\cF}_{t}$, $t\geq0$, be a filtration
of complete with respect to $\cF,P$ $\sigma$-fields
such that $ \cF _{t}\supset \tilde{\cF}_{t}$.
Our aim is to show that sometimes $u_{t}$ is
$\tilde{\cF}_{t}$-adapted even if some terms in
\eqref{2.6.4} are not 
$\tilde{\cF}_{t}$-adapted. However, the equation
is assumed to have a special structure. The result
is not surprising because in the notation, introduced
below, equation
 $$
du_{t}=
 (\Lambda^{k}_{t}u_{t}+g^{k}_{t})\,dw^{k}_{t}
$$
\begin{equation}
                                                \label{10.6.2}
+(L_{t}u_{t}+\hat{b}_{t}^{i}D_{i}u_{t}-\hat{c}_{t}u_{t}
 +D_{i}f^{i}_{t}+f^{0}_{t}+\hat{f}_{t})\,dt,
\quad t\leq\tau
\end{equation}
is written as
\begin{equation}
                                                \label{10.7.01}
du_{t}=
 (\Lambda^{k}_{t}u_{t}+g^{k}_{t})\,d\tilde{y}^{k}_{t}
+(L_{t}u_{t} 
 +D_{i}f^{i}_{t}+f^{0}_{t} )\,dt,
\quad t\leq\tau,
\end{equation}

\begin{theorem}
                                              \label{theorem 10.6.1}

Fix a number $T\in(0,\infty)$.
Assume that we are given an $\ell_{2}$-valued
process $\tilde{\sfB}_{t}$ which is $ \cF _{t}$-adapted,
jointly measurable with respect to $(\omega,t)$,
and such that $|\tilde{\sfB}_{t}|_{\ell_{2}}$
is locally  square  integrable    on $\bR_{+}$ and
$E\rho_{T}=1$, where
$$
\rho_{t} =\rho_{t}(\tilde{\sfB},dw)
=\exp(-\int_{0}^{t}
\tilde{\sfB}^{k}_{s}\,dw^{k}_{s}
-\tfrac{1}{2}\int_{0}^{t}|\tilde{\sfB}_{s}|_{\ell_{2}}^{2}\,ds). 
$$
Suppose that Assumption \ref{assumption 2.7.2} (i)
is satisfied with $\tilde{\cF}_{t}$ in place of
$\cF_{t}$ 
and
the processes
$$
\tilde{y}_{t}^{k}=w_{t}^{k}  +\int_{0}^{t}\tilde{\sfB}^{k}_{s}
\,ds
$$
are $\tilde{\cF}_{t}$ adapted.
Introduce 
$$
\bar{b}^{i}_{t}(x)=\sigma^{ik}_{t}(x)\tilde{\sfB}^{k}_{t},
\quad \bar{c}_{t}(x)=-\nu^{k}_{t}(x)\tilde{\sfB}^{k}_{t}
$$
and suppose that 
  $b+\bar{b} $ and $c+\bar{c}$ satisfy Assumption
\ref{assumption 3.11.1} with $\gamma_{b}$
from Theorem \ref{theorem 3.11.1}, for
any  $x\in
\bR^{d}$ (and $\omega$) we have $\bar{c}_{t}(x)\leq K$,
and the function
\begin{equation}
                                                         \label{9.3.1}
\int_{B_{1}}(|\bar{b}_{t}(x+y)| 
+|\bar{c}_{t}(x+y)|)\,dy
\end{equation}
is locally   integrable to the power 
$ p /(p-2)$ (locally bounded if $p=2$) on $\bR_{+}$.
Let $\tau$ be an 
$\tilde{\cF}_{t}$-stopping time such that $\tau\leq T$.

Then, for any
initial data $u_{0}\in \tr \cW^{1}_{p}(\tilde{\cF}_{0})$ and 
$f^{j},g \in\bL _{p}(\{\tilde{\cF}_{t}\},\tau)$ 
such that $\tilde{f}:=(g,\tilde{\sfB})_{\ell_{2}}
\in\bL _{p}(\{ \cF _{t}\},\tau)$,

(i)  equation \eqref{10.6.2}
has a unique solution $u$ relative to $\{ \cF _{t}\}$
in the sense of Definition \ref{definition 3.20.01},

(ii) for any $\phi\in C^{\infty}_{0}$
the process $(u_{t\wedge\tau},\phi)$
 is   $\tilde{\cF}_{t}$-adapted.
\end{theorem}

Proof. Owing to the argument
 after Theorem~\ref{theorem 1.14.1}
allowing us to introduce as large $\lambda$ as we wish,
assertion (i)   follow immediately from
  Theorem \ref{theorem 3.16.1}.

To prove (ii) we use a change of measure.  Define
$\tilde{P}(d\omega)=\rho_{T}(\omega)\,P(d\omega)$, notice that
by Girsanov's theorem the processes $\tilde{y}^{k}_{t}$, $t\leq T$,
are independent Wiener processes with respect to $\tilde{P},
\cF_{t}$.
By assumption they are $\tilde{\cF}_{t}$-adapted and since
$\tilde{\cF}_{t}\subset\cF_{t}$ the increments $\tilde{y}^{k}_{t+s}
-\tilde{y}^{k}_{t}$ are independent of $\tilde{\cF}_{t}$ if $s\geq0$.
Thus $(\tilde{y}^{k}_{t},\tilde{\cF}_{t})$ are independent Wiener
processes. Introduce $\tilde{E}$ as the expectation sign
relative to $\tilde{P}$.

After rewriting \eqref{10.6.2} in form \eqref{10.7.01}  
and applying Theorems \ref{theorem 3.11.1} and \ref{theorem 3.16.1}
we get that there exists a unique solution $\tilde{u}$ of
\eqref{10.6.2} with initial data $u_{0}$ relative to 
$\{\tilde{\cF}_{t}\}$
 in the sense of Definition
\ref{definition 3.20.01} on the new probability space, that is
with $\bL_{p}(\tau)$ and $\bW^{1}_{p}(\tau)$ replaced
with $\tilde{\,\bL}_{p}(\{\tilde{\cF}_{t}\},\tau)$ and 
$\!\tilde{\bW}^{1}_{p}(\{\tilde{\cF}_{t}\},\tau)$, respectively,
where the norms in these spaces are defined as 
$$
\tilde{E}\int_{0}^{\tau}\|u_{t}\|^{p}_{\cL_{p}}\,dt\quad
\hbox{and}\quad
\tilde{E}\int_{0}^{\tau}\|u_{t}\|^{p}_{W^{1}_{p}}\,dt
$$
raised to the power $1/p$, respectively.

Now  for $n\geq2$ we introduce $\cF_{t}$-stopping times
$$
\tau_{n}=\tau\wedge\inf\{t\geq0:\rho_{t}\leq 1/n\}
$$
and observe that 
$$
E\int_{0}^{\tau_{n}}\|\tilde{u}_{t}\|^{p}_{\cL_{p}}\,dt
\leq n E\rho_{\tau_{n}}
\int_{0}^{\tau_{n}}\|\tilde{u}_{t}\|^{p}_{\cL_{p}}\,dt=
n\tilde{E}\int_{0}^{\tau_{n}}\|\tilde{u}_{t}\|^{p}_{\cL_{p}}\,dt<\infty.
$$
Similar estimates hold if we replace $\cL_{p}$ with $W^{1}_{p}$.
By recalling that $\tilde{\cF_{t}}\subset\cF_{t}$, we conclude
that  $\tilde{u}$ is a solution of \eqref{10.6.2} 
relative to $\{\cF_{t}\}$ with $\tau_{n}$ in place of $\tau$.
By uniqueness, in the sense of distributions
$\tilde{u}_{t }I_{t\leq\tau_{n}}=
u_{t }I_{t\leq\tau_{n}}$
for almost all $(\omega,t)$, that is,
$(\tilde{u}_{t },\phi)I_{t\leq\tau_{n}}=
(u_{t },\phi)I_{t\leq\tau_{n}}$
for almost all $(\omega,t)$ for each fixed $\phi
\in C^{\infty}_{0}$. Then it follows from the integral form
of \eqref{10.6.2} that for each   $\phi
\in C^{\infty}_{0}$ with probability one
$(\tilde{u}_{t\wedge\tau_{n} },\phi) =
(u_{t \wedge\tau_{n} },\phi) $ for all $t$.
Upon  letting $n\to\infty$  we replace $\tau_{n}$
with $\tau$ and it only remains to observe that
$(\tilde{u}_{t\wedge\tau },\phi)$ is 
 $\tilde{\cF}_{t}$-measurable. The theorem is proved. 

The following is almost identical to Remark 3.5 of \cite{Kr10_1}.
\begin{remark} 

We do not use the spaces with weights.
However, there is a trivial and since very long time known
way how to use results like ours for treating equations in spaces
with weights. For instance, let $\psi_{t}(x)>0$ be a nonrandom smooth
function on $\bR^{d+1}$. Introduce, $\partial_{t}
=\partial/\partial t$,
$$
\hat{\gb}^{i}_{t}=\gb^{i}_{t}-a^{ij}_{t}D_{j}\ln\psi_{t},
\quad
\hat{b}^{i}_{t}=b^{i}_{t}-a^{ij}_{t}D_{j}\ln \psi_{t},
$$
$$
\hat{c}_{t}=c_{t}+(b^{i}_{t}+\gb^{i}_{t})D_{i}\ln\psi_{t}
-a^{ij}_{t}(D_{i}\ln\psi_{t})D_{j}\ln\psi_{t}-\partial_{t}\ln\psi_{t},
$$
$$
\hat{\nu}^{k}_{t}=\nu^{k}_{t}-\sigma^{ik}_{t}D_{i}\ln\psi_{t},
$$
$$
\hat{f}^{i}_{t}=\psi_{t} f^{i}_{t},\quad i=1,...,d,\quad
\hat{f}^{0}_{t}=f^{0}_{t}\psi_{t}-f^{i}_{t}D_{i}\psi_{t},\quad
\hat{g}^{k}_{t}=g^{k}_{t}\psi_{t}.
$$
Suppose that,
if we replace $b$, $\gb$,  $c$, and $\nu$ with $\hat{b}$,
$\hat{\gb}$, $\hat{c}$, and $\hat{\nu}$, respectively,
then
Assumptions \ref{assumption 2.7.2},  
\ref{assumption 3.11.1},  
 and \ref{assumption 3.16.1}
 are satisfied with $\gamma_{a}$
and $\gamma_{b}$ from Theorem \ref{theorem 3.11.1}.
Finally, assume that $\hat{f}^{j},\hat{g}
\in\bL_{2}(\tau)$ and $u_{0}\psi_{0}\in\tr\cW^{1}_{p}$. 
Then it turns out that for $\lambda\geq\lambda_{0}$
($\lambda_{0}$ is taken from Theorem \ref{theorem 3.11.1})
equation
\eqref{2.6.4} has a unique solution $u$ such that
$u\psi\in\bW^{1}_{p}(\tau)$.

This fact is almost trivial since $u$ satisfies 
\eqref{2.6.4} if and only if $v:=u\psi$
satisfies the version of  
\eqref{2.6.4} which is obtained
as the result of the replacements described above
and also the replacement of $f^{j},g$ with $\hat{f}^{j},\hat{g}$,
respectively.
In addition,
the natural estimate of the $\bW^{1}_{p}(\tau)$-norm of $v$ gives
an estimate of $u$ in an appropriate space with weights.

As a specification of the above, in the setting of Remark
\ref{remark 1.15.1} take a $T\in(0,\infty)$, set $\tau=T$,
and for $\theta\in(0,\infty)$  introduce
$$
\ln \psi_{t}(x)=-\theta e^{\theta^{2}(t-T)}\sqrt{1+|x|^{2}}.
$$
Obviously, $D_{i}\ln\psi$ are bounded for $t\leq T$.
Furthermore, it is not hard to see that if $\theta$ is large enough,
then $\hat{c}_{t}\geq0$ for $t\leq T$. Also,
if $|x-y|\leq1$, then owing to the fact that $|D_{ij}\ln\psi_{t}(x)|
\leq N(1+|x|)^{-1}$ for $t\leq T$, where $N$ is a constant, we have
$$
|b^{i}_{t}(x)D_{i}\ln\psi_{t}(x)-b^{i}_{t}(y)D_{i}\ln\psi_{t}(y)|
$$
$$
\leq
|(b^{i}_{t}(x)-b^{i}_{t}(y))D_{i}\ln\psi_{t}(x)|
+N(1+|x|)|D\ln\psi_{t}(x)- D\ln\psi_{t}(y)|
$$
$$
\leq K|D\ln\psi_{t}(x)|+N
$$
for $t\leq T$.
Estimates similar to this one show that $\hat{b}$,
$\hat{\gb}$, and $\hat{c}$ satisfy
Assumption \ref{assumption 3.11.1}
 for $t\leq T$. By what is said
in the beginning of the current remark, if $u_{0}\psi_{0}\in
\tr\cW^{1}_{p}$ (for instance,  $u_{0}(x)=x^{1}$),
then \eqref{2.6.4} has a unique solution $u$ such that
$u\psi\in \bW^{1}_{p}(T)$. Since $D\ln\psi$ is bounded,
the inclusion $u\psi\in \bW^{1}_{p}(T)$ is equivalent
to $u\psi\in \bL_{p}(T)$, $\psi Du\in \bL_{p}(T)$.

To the best of the author's knowledge 
 even in this special case the result 
in this generality was not known before.

\end{remark}

\mysection{Proof of Theorems \protect\ref{theorem 3.11.1}
and \protect\ref{theorem 1.12.1}}
                                              \label{section 1.1.1}

In this section we suppose that Assumptions
\ref{assumption 2.7.2}, \ref{assumption 2.4.1},
\ref{assumption 2.4.01}, and
\ref{assumption 3.11.1} are satisfied with some $\gamma_{a},
\gamma_{b}\in(0,1]$ and 
 start by showing that the requirement (i)
of Definition \ref{definition 3.20.01} is automatically satisfied
for any $u\in\bW^{1}_{p}(\tau)$. Take  a  nonnegative
$ \xi\in C^{\infty}_{0}(B_{\rho_{1}})$
  with unit integral and 
define
$$
\bar{b}_{s}(x)=\int_{B_{\rho_{1}}}\xi(y) b_{s}(x-y) \,dy,\quad
\bar{\gb}_{s}(x) =\int_{B_{\rho_{1}}}\xi(y) \gb_{s}(x-y) \,dy,
$$
\begin{equation}
                                                      \label{6.28.3}
\bar{c}_{s}(x)=\int_{B_{\rho_{1}}}\xi(y) c_{s}(x-y) \,dy.
\end{equation}
We may assume that $|\xi|\leq N(d)\rho_{1}^{-d}$.
\begin{remark}
                                               \label{remark 1.12.1}
By Corollary 5.4 of \cite{Kr10_2} , for 
$x_{0}\in\bR^{d}$, $v\in\cL_{p}$,
$\phi\in W^{1}_{p'}$, and $u\in W^{1}_{p}$ we have
$$
 (|b _{s}-\bar{b} _{s}(x_{0})|I_{B_{\rho_{1}}(x_{0})}
v,|\phi|) 
\leq N \| v\| _{\cL_{p}}
\|\phi\|_{W^{1}_{p'}} ,
$$  
\begin{equation}
                                                   \label{1.12.1}
\|I_{B_{\rho_{1}}(x_{0})}
|\gb_{s}-\bar{\gb}_{s}(x_{0})|\,u \|_{\cL_{p}}
+\|I_{B_{\rho_{1}}(x_{0})}|c_{s}-\bar{c}_{s}(x_{0})|\,u\|_{\cL_{p}}
\leq N\|u\|_{W^{1}_{p}},
\end{equation}
where $N=N(d,p,\rho_{1},K) $. In particular,
\begin{equation}
                                                   \label{1.12.2} 
(|b _{s}|I_{B_{\rho_{1}}(x_{0})}
v,|\phi|) 
\leq (N+|\bar{b} _{s}(x_{0})|)
 \| v\| _{\cL_{p}}
\|\phi\|_{W^{1}_{p'}} ,
\end{equation}
the latter implying that $|b _{s}|I_{B_{\rho_{1}}(x_{0})}
v\in H^{-1}_{p}$. It is also seen that
if $u\in\bW^{1}_{p}(\tau)$ and $|\bar{b} _{s}(x_{0})|$
is a {\em bounded\/} function on $\opar0,\tau\cbrk$,
then 
$$
I_{B_{\rho_{1}}}(x_{0})b^{i} D_{i}u 
\in\cL_{p}(\opar0,\tau\cbrk,\cP,H^{-1}_{p}).
$$
Similarly,
\begin{equation}
                                                   \label{1.13.1}
\|I_{B_{\rho_{1}}(x_{0})}
|\gb_{s}|\,u \|_{\cL_{p}}
+\|I_{B_{\rho_{1}}(x_{0})}|c_{s}|\,u\|_{\cL_{p}}
\leq (N+|\bar{\gb}_{s}(x_{0})|+
|\bar{c}_{s}(x_{0})|)\|u\|_{W^{1}_{p}}.
\end{equation}

By the way, Remark \ref{remark 1.12.3} now shows that
under the conditions of Theorem \ref{theorem 1.12.1}
for any solution $u$ of \eqref{2.6.4} and $\phi\in
C^{\infty}_{0}$ with support lying in a ball of radius
$\rho_{1}$ we have $u\phi\in\cW^{1}_{p}(\tau)$. Of course,
the restriction on the size of support of $\phi$ is easily
removed and this proves assertion (i) of
Theorem \ref{theorem 1.12.1}.

\end{remark}

\begin{lemma}
                                  \label{lemma 6.27.1}
Let $R\in(0,\infty)$. Then there exists a sequence
of bounded stopping times $\tau_{n}\to\infty$ such that for any
$\omega\in\Omega$,
  $u\in\cL_{p}((0,\tau),W^{1}_{p})$, and  
 $\phi\in C^{\infty}_{0}(B_{R})$
$$
\int_{0}^{\tau_{n}\wedge\tau}
 (|(b^{i}_{s}D_{i}u_{s},\phi)|+
 |(\gb^{i}_{s} u_{s},D_{i}\phi)|+
|(c_{s}u_{s},\phi)|)\,ds
$$
\begin{equation}
                                             \label{8.12.3}
\leq  n\|u\|_{\cL_{p}((0,\tau),W^{1}_{p})}
\|\phi\|_{W^{1}_{p'}},
\end{equation}
so that requirement (i) in Definition \ref{definition 3.20.01}
can be dropped.  

\end{lemma}
 Proof. By having in mind partitions of unity we convince
ourselves that it suffices to prove \eqref{8.12.3}
under the assumption that $\phi$ has support in a ball
$B_{\rho_{1}}(x_{0})$. 
Observe that  by \eqref{1.13.1} and H\"older's inequality
\begin{equation}
                                            \label{1.15.4}
|(\gb^{i}_{s} u_{s},D_{i}\phi)|+|(c_{s}u_{s},\phi)|\leq
N(1+|\bar{\gb}_{s}(x_{0})|+
|\bar{c}_{s}(x_{0})|)\|u_{s}\|_{W^{1}_{p}}\|\phi\|_{W^{1}_{p'}}.
\end{equation}
It follows again by H\"older's inequality that
$$
\int_{0}^{t\wedge\tau}(
 |(\gb^{i}_{s} u_{s},D_{i}\phi)|+|(c_{s}u_{s},\phi)|)\,ds
\leq N \chi_{t}\|u\|_{\cL_{p}((0,\tau),W^{1}_{p})}
\|\phi\|_{W^{1}_{p'}},
$$
where
$$
\chi_{t}= t^{1/p'}
+\big(\int_{0}^{t}|\bar{\gb} _{s}(x_{0})|^{p'}ds\big)^{1/p'}
+\big(\int_{0}^{t}|\bar{c} _{s}(x_{0})|^{p'}\,
ds\big)^{1/p'} .
$$
After that, in what concerns $\gb$ and $c$, it only remains to recall
Assumption \ref{assumption 2.7.2} (iii). 
Similarly the integral of $ |(b^{i}_{s}D_{i}u_{s},\phi)|$
is estimated by using \eqref{1.12.2} and the lemma is proved.
\begin{remark}
                                           \label{remark 9.1.1}
Estimates \eqref{1.12.2} and \eqref{1.13.1}
show that for any $u\in\bW^{1}_{p}$ for almost all
$(\omega,s)$ the functions $b^{i}_{s}D_{i}u_{s}$,
$D_{i}(\gb^{i}_{s}u_{s})$, and $c_{s}u_{s}$
 are distributions on $\bR^{d}$.
\end{remark}

Since bounded linear operators
are continuous we obtain the following.
 \begin{corollary}
                                    \label{corollary 3.23.1}
Let $R,\tau_{n},\phi$ be as in Lemma \ref{lemma 6.27.1}. Then
 the operators
$$
u_{t}\to \int_{0}^{t\wedge\tau_{n}}
(b^{i}_{s}D_{i}u_{s},\phi)\,ds,\quad
u_{t}\to \int_{0}^{t\wedge\tau_{n}}
(\gb^{i}_{s}u_{s},D_{i}\phi)\,ds, 
$$
$$
u_{t}\to \int_{0}^{t\wedge\tau_{n}}
(c_{s}u_{s}, \phi)\,ds
$$
are continuous as operators from $\bW^{1}_{p}(\tau)$ to
$\cL_{p}(\opar0,\tau_{n}\cbrk)$ for any $n$.
\end{corollary}  
This result will be used in Section
\ref{section 6.9.5}.

 Now we   prove Theorem \ref{theorem 3.11.1}
in a particular case.

\begin{lemma}
                                        \label{lemma 3.11.1}
Let  $\gb^{i}$, $b^{i}$, and $c$ 
be independent of $x$ and let  $u_{0}=0$. 
Then the assertion of Theorem \ref{theorem 3.11.1}
holds, naturally, with $\lambda_{0}=\lambda_{0}(d,\delta, p,
\rho_{0},\kappa )$ (independent of $\rho_{1}$).

\end{lemma}

Proof.  First let  $c\equiv0$.
We want to use the It\^o-Wentzell formula to get rid
of the first-order terms. Observe that \eqref{2.6.4} reads as
$$
du_{t}=(\Lambda_{t}^{k}u_{t}+g^{k}_{t})\,dw^{k}_{t}
$$
\begin{equation}
                                                \label{6.28.1}
 +\big[D_{i}(a^{ij}_{t}D_{j}u_{t}+(\gb^{i}_{t} + 
b^{i}_{t})u_{t}+f^{i}_{t})+f^{0}_{t}- \lambda u_{t}\big]\,dt,
\quad t\leq \tau.
\end{equation}
 
Recall that from the start
(see Definition \ref{definition 3.20.01})
 it is assumed that $u \in  
\bW^{ 1}_{p}(\tau)$. Then
one can find a predictable set $A\subset\opar0,\tau\cbrk$
of full measure such that $I_{A}f^{j}$, $j=0,1,...,d$,
$I_{A} g$, and
$I_{A}D_{i}u$, $i=1,...,d$,
 are well defined as $\cL_{p}$-valued predictable
functions satisfying
$$
\int_{0}^{\tau}I_{A}\big(\sum_{j=0}^{d}\|f^{j}_{t}\|^{p}_{\cL_{p}}
+ \|g_{t}\|^{p}_{\cL_{p}}+ \|Du_{t}\|^{p}_{\cL_{p}}\big)\,dt<\infty.
$$
Replacing $f^{j}$, $g$, and $D_{i}u$ in \eqref{6.28.1} with
$I_{A}f^{j}$, $I_{A}g $,
 and $I_{A}D_{i}u$, respectively, will not affect
\eqref{6.28.1}. Similarly one can treat the term 
$h_{t}= (\gb^{i}_{t}+b^{i}_{t})u_{t}$
for which
$$
\int_{0}^{T\wedge\tau}\|h_{t}\|_{\cL_{p}}\,dt<\infty
$$
(a.s.) for each   $T\in\bR_{+}$, owing to Assumption 
\ref{assumption 2.7.2} and 
the fact that $u\in\bL_{p}( \tau)$.

After these replacements all terms on the right in \eqref{6.28.1}
will be of   class $\frD^{1}  $ and $\frD^{2} $ as appropriate
  since  
 $a$ and $\sigma$ are bounded (see the definition of 
$\frD^{1} $ and $\frD^{2} $
in \cite{Kr09_4}). 
This allows us to apply Theorem 1.1 of \cite{Kr09_4}
and for
$$
B_{t}^{i}=\int_{0}^{t}
(\gb^{i}_{s}+b^{i}_{s})\,ds,\quad 
 \hat{u}_{t}(x)=u_{t}(x-B_{t})
$$
obtain that
\begin{equation}
                                               \label{6.28.8}
d\hat{u}_{t}=\big[D_{i}(\hat{a}^{ij}_{t}D_{j}\hat{u}_{t} )
 - \lambda \hat{u}_{t}+D_{i}\hat{f}^{i}_{t}+\hat{f}^{0}_{t}
\big]\,dt+(\hat{\Lambda}^{k}_{t}\hat{u}_{t}+\hat{g}^{k}_{t})
\,dw^{k}_{t},
\end{equation}
where $\hat{\Lambda}^{k}_{t}=\hat{\sigma}^{ik}_{t}D_{i}+
\hat{\nu}^{k}_{t}$ and
$$
(\hat{a}^{ij}_{t},\hat{\sigma}^{ik}_{t},
\hat{\nu}^{k}_{t}
\hat{f}^{j}_{t},\hat{g}^{k}_{t})(x)
=(a^{ij}_{t}, \sigma^{ik}_{t},\nu^{k}_{t}, f^{j}_{t},
g^{k}_{t})(x-B_{t}).
$$

Obviously, $\hat{u}$ is in $\bW^{1}_{p}(\tau)$ and its norm
coincides with that of $u$. 
Equation \eqref{6.28.8}  
shows that $\hat{u}\in\cW^{1}_{p}(\tau)$.

Next observe that owing to    \eqref{2.9.5},
for any $\omega\in\Omega,\rho\in(0,\rho_{0}],
t\geq0$, and $i,j=1,...,d$ we have 
$$                                           
\rho^{-2d-2}\int_{t}^{t+\rho^{2}}\bigg(
 \sup_{x\in \bR^{d}}\int_{B_{\rho}(x)}\int_{B_{\rho}(x)}
|\hat{a}^{ij}_{s}(  y)-\hat{a}^{ij}_{s}( z)|\,dydz\bigg)
 \,ds\leq\gamma_{a}, 
$$
which in terms of \cite{Kr09} implies that
the couple $(\hat{a},\hat{\sigma})$ is $(\varepsilon,
\varepsilon )$-regular
at any point of $\bR_{+}\times\bR^{d}$
for any $\varepsilon\in(0,\rho_{0}]$. 
Then owing to our Assumptions \ref{assumption 2.7.2}
(ii) and \ref{assumption 2.4.1}
one can choose  $\varepsilon=\varepsilon(\delta,\kappa)
\in(0,\rho_{0}]$ so that Assumption 2.2 of \cite{Kr09}
is satisfied.

By Theorem 2.2
of \cite{Kr09} if 
Assumption \ref{assumption 2.4.01}
is satisfied with $\gamma_{a}=\gamma 
_{a}(d,\delta,p)>0$, specified in its proof,
and if
$\lambda\geq\lambda_{0}(d,\delta,p,\kappa,
\rho_{0})\geq1$,
then
$$
\lambda\|\hat u\|^{2}_{\bL _{p}(\tau)}+\|D\hat u\|^{2}_{\bL _{p}(\tau)}
\leq N\big(\sum_{i=1}^{d}\|\hat f^{i}\|^{2}_{\bL _{p}(\tau)}
+\|\hat g\|^{2}_{\bL _{p}(\tau)} 
+ \lambda^{-1}\|\hat f^{0}\|^{2}_{\bL _{p}(\tau)}\big) ,
$$
where $N=N(d,\delta,p,\kappa,
\rho_{0})$. This coincides with \eqref{3.11.2}
and proves the lemma in case $c\equiv0$.

In the general case observe that 
owing to Assumption \ref{assumption 2.7.2} (iii)
there exists a sequence of stopping times $\tau_{n}\uparrow\tau$
such that 
$$
\int_{0}^{\tau_{n}}c_{s}\,ds\leq n.
$$
Clearly, if we can prove \eqref{3.11.2} with $\tau_{n}$
in place of $\tau$, then by passing to the limit we will get
\eqref{3.11.2} as is. Therefore, without losing generality
we assume that
$$
\sup_{\Omega}\int_{0}^{\infty}c_{s}\,ds<\infty.
$$

Then introduce
$$
\xi_{t}=\exp(\int_{0}^{t}c_{s}\,ds).
$$
By the above argument    we have $\bar{u}:=
\xi u\in \bW^{1}_{p}(\tau)$ and
$$
d\bar{u}_{t}= \big[D_{i}(a^{ij}_{t}D_{j}\bar{u}_{t}+[\gb^{i}_{t} + 
b^{i}_{t}]\bar{u}_{t}+\xi_{t}f^{i}_{t})+\xi_{t}f^{0}_{t}- \lambda 
\bar{u}_{t}
\big]\,dt+(\Lambda^{k}_{t}\bar{u}_{t}
+\xi_{t}g^{k}_{t})\,dw^{k}_{t},
\quad t\leq \tau.  
$$
By the above result for any stopping time   $\tau'\leq \tau$
$$
\lambda^{p/2}\|\xi u\|^{p}_{\bL _{p}(\tau')}+
\|\xi Du\|^{p}_{\bL _{p}(\tau')}=
\lambda^{p/2}\|\bar{u}\|^{p}_{\bL _{p}(\tau')}+
\|D\bar{u}\|^{p}_{\bL _{p}(\tau')}
$$ 
\begin{equation}
                                       \label{3.11.02}
\leq N\big(\sum_{i=1}^{d}\|\xi f^{i}\|^{p}_{\bL _{p}(\tau')}
+\|\xi g\|^{p}_{\bL _{p}(\tau')} 
+ \lambda^{-p/2}\|\xi f^{0}\|^{p}_{\bL _{p}(\tau')}\big).
\end{equation}

If needed, one can enlarge the original probability
space in such a way that there will exist an
exponentially distributed, with parameter one,
random variable $\eta$ independent of $\{\cF_{t},t\geq0\}$.
 We assume that the
enlargement is not needed and
define
$$
\phi_{t}=p\int_{0}^{t}c_{s}\,ds,\quad
\psi_{s}=\tau\wedge\inf\{t\geq 0: \phi_{t}\geq s\},
\quad\tau' =\psi_{\eta}.
$$
Notice that
$$
\{\omega:\psi_{s}> t\}=\{\omega:\tau>t, \phi_{t}< s\}.
$$
Hence  
$$
\{\omega:\tau'> t\}=\{\omega:\tau>t,  \phi_{t}< \eta\}.
$$
It follows that $\tau'$ is a stopping time
with respect to $\cF_{t}\vee\sigma(\eta)$. Furthermore,
for any nonnegative  predictable (relative to
the original filtration $\cF_{t}$) process $h_{t}$
we have
$$
E\int_{0}^{\tau'}h_{t}\,dt=\int_{0}^{\infty}
Eh_{t}E\{I_{\tau'>t}\mid \cF_{t}\}\,dt
$$
$$
=\int_{0}^{\infty}
Eh_{t} I_{\tau >t}e^{- \phi_{t}}\,dt=E\int_{0}^{\tau}h_{t}
\xi^{-p}_{t}\,dt.
$$

This  and \eqref{3.11.02}
immediately lead to \eqref{3.11.2}
and the lemma is proved.

  To proceed further take $\bar{\gb},
\bar{b}$, and $\bar{c}$
from \eqref{6.28.3}.
 From Lemma 4.2 of \cite{Kr09_1} and Assumption
\ref{assumption 3.11.1} it follows  that, for 
$h_{ t}=\bar{\gb}_{ t},\bar{b}_{  t},
\bar{c}_{ t}$, it holds that  
 $|D^{n}h_{ t}|\leq M_{n} $, where
$M_{n}=M_{n}
(n, d,\rho_{1},K)\geq1$ and $ D^{n}h_{ t }$
is any derivative of $h_{ t}$ of order $n\geq1$
with respect to $x$. By Corollary 4.3 of \cite{Kr09_1}
we have $|h_{ t}( x )|\leq K(t)(1+|x|)$,
where for each $\omega$ the function $K(t)=K(\omega,t)$
is locally 
 integrable with respect to $t$ on $\bR_{+}$.
Owing to these properties
the equation
\begin{equation}
                                             \label{2.8.1}
x_{t}=x_{0}-\int_{t_{0}}^{t}(\bar{\gb}_{ s}+\bar{b}_{  s})
( x_{s})\,ds,\quad t \geq t_{0} ,
\end{equation}
 for any ($\omega$ and)  
$ (t_{0},x_{0}) \in \bR^{d+1 }_{+} $  has a unique solution 
$x_{t}=x_{t_{0},x_{0},t} $.
Obviously, the process $x_{t_{0},x_{0},t}$, $t\geq t_{0}$,
 is $\cF_{t}$-adapted.

Next, for $i=1,2$ set $\chi^{(i)}(x)$ to be the indicator function
of $B_{\rho_{1}/i}$ and introduce
$$
\chi^{(i)}_{t_{0},x_{0},t}(x)=\chi^{(i)}(x-x_{t_{0},x_{0},t})
I_{t\geq t_{0}}.
$$

By using the above results
and reproducing  the proofs of Lemma 5.5 of \cite{Kr10_1},
where $p=2$ and SPDEs are treated,
and Lemma 5.8 of~\cite{Kr10_2}, where $p$ is general
but only PDEs are considered, we easily obtain the following.
\begin{lemma}
                                     \label{lemma 3.14.1}
Suppose that
Assumption \ref{assumption 2.4.01}
is satisfied with   $\gamma_{a}=\gamma_{a}(d,\delta,p)$
taken from Lemma \ref{lemma 3.11.1}.
Assume that   we are given a function
$u  $ which is a solution of \eqref{2.6.4} with  
   some
$f^{j},g \in\bL_{p}(\tau)$, and $\lambda\geq\lambda_{0}=
\lambda_{0}(d,\delta,p,  
\rho_{0},\kappa )$, where $\lambda_{0}(d,\delta, p,
\rho_{0},\kappa )$ is taken from Lemma~\ref{lemma 3.11.1}.
Take $ (t_{0},x_{0}) \in 
\bR^{d+1}_{+}$   
 and assume that $u_{t}=0$ if $t\leq t_{0}\wedge \tau$.
 Then  
$$
 \lambda^{p/2} \|\chi^{(2)}_{t_{0},x_{0}}u \|^{p}_{\bL _{p}(\tau)}+
\|\chi^{(2)}_{t_{0},x_{0}}Du\|^{p}_{\bL _{p}(\tau)} 
$$
$$
\leq N\big(\sum_{i=1}^{d}\|\chi^{(1)}_{t_{0},x_{0}}
f^{i}\|^{p}_{\bL _{p}(\tau)}+\|
\chi^{(1)}_{t_{0},x_{0}}g\|^{p}_{\bL_{p}(\tau)}\big)
+N\lambda^{-p/2}\|\chi^{(1)}_{t_{0},x_{0}}f^{0}\|^{p}_{\bL _{p}(\tau)}
$$
$$
+N\gamma_{b}^{p/q} \| 
\chi^{(1)}_{t_{0},x_{0}} Du
\|_{\bL_{p}(\tau)}^{p}+
 N^{*} \lambda^{-p/2}\| 
\chi^{(1)}_{t_{0},x_{0}} Du
\|_{\bL_{p}(\tau)}^{p}
$$
\begin{equation}
                                       \label{3.14.2}
+ N^{*}  \|
\chi^{(1)}_{t_{0},x_{0}}  u
\|_{\bL_{p}(\tau)}^{p}
+N^{*}\lambda^{-p/2}\sum_{i=1}^{d}\|\chi^{(1)}_{t_{0},x_{0}}
f^{i}\|^{p}_{\bL _{p}(\tau)},
\end{equation}
where  
   $N$ is a constant depending only
on $d,\delta$, $p$, $\rho_{0}$,   and
$\kappa$  and   $N^{*}$   depends only
on the same objects, $\gamma_{b}$, $\rho_{1}$,  and $K$.
\end{lemma}

Upon integrating through equation \eqref{3.14.2}
with respect to $x_{0}$ and repeating the arguments
in the proofs of Lemma 5.6 of \cite{Kr10_1}
or Lemma 5.9 of \cite{Kr10_2} we obtain the following
result in which $M_{1}( d,\rho_{1},K)$ is the constant introduced
before Lemma \ref{lemma 3.14.1}.

\begin{lemma}
                                                \label{lemma 3.14.3}

Suppose that
Assumption \ref{assumption 2.4.01}
is satisfied with   $\gamma_{a}=\gamma_{a}(d,\delta,p) $
taken from Lemma \ref{lemma 3.11.1}.
Assume that   we are given a function
$u  $ which is a solution of \eqref{2.6.4} with  
   some
$f^{j},g \in\bL_{p}(\tau)$, and $\lambda\geq\lambda_{0}=
\lambda_{0}(d,\delta,p,  
\rho_{0},\kappa )$, where $\lambda_{0}(d,\delta, p,
\rho_{0},\kappa )$ is taken from Lemma~\ref{lemma 3.11.1}.
 Take an  $s_{0}\in\bR_{+}$ 
 and assume that $u_{t}=0$ if $t\leq s_{0}\wedge\tau$.
 Then for   $I_{s_{0}}:=
I_{(s_{0},t_{0})}$, where $t_{0}= s_{0}+M_{1}^{-1} $,
we have
$$
\lambda^{p/2}\| I_{s_{0}}u\|^{p}_{\bL _{p}(\tau) }+
\|I_{s_{0}} Du\|^{p}_{\bL _{p}(\tau) }
\leq N\big( \sum_{i=1}^{d}\| 
I_{s_{0}}f^{i}\|^{p}_{\bL _{p}(\tau) }+\| 
I_{s_{0}}g\|^{p}_{\bL _{p}(\tau) }\big)
$$
$$
+N\lambda^{-p/2 }\|I_{s_{0}} f^{0}\|^{p}_{\bL _{p} (\tau)}
 +N\gamma_{b} ^{p/q} \| I_{s_{0}} Du\|_{\bL_{p}(\tau) }^{p}
$$
\begin{equation}
                                       \label{3.14.5} 
+
 N^{*} \lambda^{-p/2}\| I_{s_{0}}  Du\|_{\bL_{p}(\tau) }^{p}
+ N^{*} \| I_{s_{0}} u\|_{\bL_{p}(\tau) }^{p}
+N^{*}\lambda^{-p/2}\sum_{i=1}^{d}\| 
I_{s_{0}}f^{i}\|^{p}_{\bL _{p}(\tau) },
\end{equation}
where  
   $N$ is a constant depending only
on $d,\delta$, $p$, $\rho_{0}$,   and
$\kappa$  and   $N^{*}$   depends only
on the same objects, $\gamma_{b}$, $\rho_{1}$, and $K$.
\end{lemma}

{\bf  Proof of Theorem \ref{theorem 3.11.1}}. First we show 
how to choose
an appropriate $\gamma_{b}=\gamma_{b}(d,\delta,p,\rho_{0},\kappa)$.
Call $N_{0}$ the constant factor of $\gamma_{b}^{p/q} \| 
I_{s_{0}} Du\|_{\bL_{p}(\tau) }^{p}$ in \eqref{3.14.5} and choose 
a $\gamma_{b}\in(0,1]$ in such a way that $N_{0}
\gamma^{p/q}_{b}\leq 1/2$. Then under the assumptions of Lemma
\ref{lemma 3.14.3} we have
$$
\lambda^{p/2}\| I_{s_{0}}u\|^{p}_{\bL _{p}(\tau) }+
\|I_{s_{0}} Du\|^{p}_{\bL _{p}(\tau) }
\leq N \big(\sum_{i=1}^{d}\| 
I_{s_{0}}f^{i}\|^{p}_{\bL _{p}(\tau) }+\| 
I_{s_{0}}g\|^{p}_{\bL _{p}(\tau) }\big)
$$
$$
+N\lambda^{-p/2 }\|I_{s_{0}} f^{0}\|^{p}_{\bL _{p} (\tau)}
+ N^{*} \lambda^{-p/2}\| I_{s_{0}}  Du\|_{\bL_{p}(\tau) }^{p}
$$
\begin{equation}
                                       \label{8.28.1} 
+ N^{*} \| I_{s_{0}} u\|_{\bL_{p}(\tau) }^{p}
+N^{*}\lambda^{-p/2}\sum_{i=1}^{d}\| 
I_{s_{0}}f^{i}\|^{p}_{\bL _{p}(\tau) }.
\end{equation}

To proceed further assume that
\begin{equation}
                                       \label{8.28.3} 
u_{0}=0.
\end{equation}

After $\gamma_{b}$ has been fixed we recall that $M_{1}=M_{1}
( d,\rho_{1},K)$
 and we take a $\zeta\in C^{\infty}_{0}(\bR)$  
with support in $(0,M_{1}^{-1})$
such that
\begin{equation}
                                       \label{3.15.1} 
\int_{-\infty}^{\infty}\zeta^{p}(t)\,dt=1.
\end{equation}
 For $s\in\bR $ define
$\zeta^{s}_{t}=\zeta(t-s)$, 
$u^{s}_{t}( x)=u_{t}(x)\zeta^{s}_{t}$. 
Obviously, $u^{s}_{t}=0$ if $0\leq t\leq s_{+}\wedge\tau$.
Therefore, we can apply 
\eqref{8.28.1}  to $u^{s}_{t}$  by taking $s_{0}=s_{+}$ and
 observing that
$$
du^{s}_{t}=(L_{t}u^{s}_{t}-\lambda u^{s}_{t}
+D_{i}(\zeta^{s}_{t}f^{i}_{t})+\zeta^{s}_{t}f^{0}_{t}+u_{t}(\zeta^{s}_{t})')\,dt
+(\Lambda^{k}_{t}u^{s}_{t}+\zeta^{s}_{t}g^{k}_{t})\,dw^{k}_{t},
\quad t\leq\tau.
$$
We also use the fact that for $t\geq0$, as is easy to see,  
 $I_{s_{+}}(t)\zeta^{s}_{t}=\zeta^{s}_{t}$.
Then for  and $\lambda\geq\lambda_{0}=
\lambda_{0}(d,\delta,p,  
\rho_{0},\kappa )$, where $\lambda_{0}(d,\delta, p,
\rho_{0},\kappa )$ is taken from Lemma~\ref{lemma 3.11.1},
we obtain
$$
\lambda^{p/2}\| \zeta^{s} u\|^{p}_{\bL _{p}(\tau) }+
\|\zeta^{s}  Du\|^{p}_{\bL _{p}(\tau) }
\leq N \big(\sum_{i=1}^{d}\| 
\zeta^{s} f^{i}\|^{p}_{\bL _{p}(\tau) }
+\| 
\zeta^{s}g\|^{p}_{\bL _{p}(\tau) }\big)
$$
$$
+N\lambda^{-p/2 }(\|\zeta^{s}  f^{0}\|^{p}_{\bL _{p} (\tau)}
+\|(\zeta^{s})' u\|^{p}_{\bL _{p} (\tau)})
$$
\begin{equation}
                                       \label{8.28.2} 
+ N^{*} \lambda^{-p/2}\| \zeta^{s}  Du\|_{\bL_{p}(\tau) }^{p}
+ N^{*} \| \zeta^{s} u\|_{\bL_{p}(\tau) }^{p}
+N^{*}\lambda^{-p/2}\sum_{i=1}^{d}\| 
\zeta^{s} f^{i}\|^{p}_{\bL _{p}(\tau) }.
\end{equation}
We integrate through this relation with respect to
$s\in\bR $, use \eqref{3.15.1}
and 
$$
\int_{-\infty}^{\infty}|(\zeta^{s}_{t})'|^{p}\,ds=
\int_{-\infty}^{\infty}| \zeta'(t) |^{p}\,dt=N^{*}.
$$

Then we conclude
$$
\lambda^{p/2}\| u\|^{p}_{\bL _{p}(\tau) }+
\|  Du\|^{p}_{\bL _{p}(\tau) }
\leq N_{1} \big(\sum_{i=1}^{d}\| 
  f^{i}\|^{p}_{\bL _{p}(\tau) }
+\| g\|^{p}_{\bL _{p}(\tau) }\big)
$$
$$
+ N_{1}\lambda^{-p/2 } \|   f^{0}\|^{p}_{\bL _{p} (\tau)}
+ N^{*} _{1} \lambda^{-p/2}\|   Du\|_{\bL_{p}(\tau) }^{p}
$$
$$
+ N^{*} _{1} \|  u\|_{\bL_{p}(\tau) }^{p}
+N^{*} _{1}\lambda^{-p/2}\sum_{i=1}^{d}\| 
 f^{i}\|^{p}_{\bL _{p}(\tau) }.
$$
Without losing generality we assume that $N_{1}\geq1$
and  we show how to choose $\lambda_{0}=\lambda_{0}(
d,\delta,p,\rho_{0},\rho_{1},\kappa,K)\geq1$. 
Above we assumed that $\lambda\geq\lambda_{0}(d,\delta,p,\rho_{0}
,\kappa)$,
where $\lambda_{0}(d,\delta,p,\rho_{0},\kappa)$ is taken from
Lemma \ref{lemma 3.11.1}. Therefore, we take
$$
\lambda_{0}=\lambda_{0}(
d,\delta,p,\rho_{0},\rho_{1},\kappa,K)\geq
\lambda_{0}(d,\delta,p,\rho_{0},\kappa)
$$
  such that
  $\lambda_{0}^{p/2}\geq 2N^{*}_{1} $
(recall that $N_{1}^{*}=N_{1}^{*}
(d,\delta,p,\rho_{0},\rho_{1},\kappa,K)$). Then we obviously come
to \eqref{3.11.2} (with $u_{0}=0$). 

A standard method to remove assumption \eqref{8.28.3}
by subtracting from $u$ the solution of the heat equation
$dv_{t}=(\Delta v_{t}-v_{t})\,dt$ with initial data $u_{0}$
does not work because it leads to subtracting
the terms $D_{i}(\gb^{i}v)+b^{i}D_{i}v$, which one should include
into the free terms  $D_{i}f^{i}+f^{0}$ 
in the equation. Generally, this is 
impossible because we only know that $D_{i}v\in\bL_{p}(\tau)$
and if we multiply $D_{i}v$ by an arbitrary 
function of $x$ with linear growth,
the inclusion may fail.

Therefore, we use a different method.
The idea is to shift
all data along the time axis by 1, consider our equations
on $\opar1,\hat \tau \cbrk$, where $\hat \tau =1+\tau$, and 
supplement  this equation with
an equation for $t\in[0,1]$ with zero initial data and such that
the value of its solution at time 1 would coincide with $u_{0}$.
Then the two equations combined would give an equation on   
$\opar0,\hat \tau \cbrk$ with zero initial condition, which would allow
us to apply the above result.

Formally, we need to have   Wiener processes on $[0,\infty)$
and after shifting they will be defined only on $[1,\infty)$
(and satisfy $w^{k}_{1}=0$). Therefore, we augment if needed
our probability space in such a  way that we may assume that
there are Wiener processes $\bar{w}^{1}_{t},
\bar{w}^{2}_{t},...$, $t\geq0$, independent of $\{\cF_{s},s\geq0\}$.
Then define $\cF^{\bar w}_{t}$ as the completion
of $\sigma(\bar{w}_{s}:s\leq t)$,
$$
\hat\cF _{t}=\cF_{0}\vee\cF^{\bar w}_{t},\quad t\in[0,1],
\quad \hat\cF _{t}=\cF_{t-1}\vee\cF^{\bar w}_{1},\quad t\geq1,
$$
$$
\hat w^{ k}_{t}=\bar{w}^{k}_{t},\quad t\in[0,1],\quad
\hat w^{ k}_{t}=\bar{w}^{k}_{1}+w^{k}_{t-1}\quad t\geq1,
\quad
\hat\tau=1+\tau,
$$
and for $t\geq1$ define the coefficients and the free terms
by following the example $\hat{a}^{ij}_{t}=a^{ij}_{t-1}$.

Next, take the function $v$ from Definition \ref{definition 8.10.1}
and for $t\in[0,1]$  set 
$$
\hat{a}^{ij}_{t}=\delta^{ij},\quad \hat f^{i}_{t}=-2tD_{i}v_{1-t},
\quad\hat f^{0}_{t}=(1+t+\lambda t)v_{1-t},
$$
where $\lambda\geq\lambda_{0}$ with $\lambda_{0}$ determined
in the first part of the proof. We define all other coefficients
with hats and
the free terms $\hat{g}^{k}_{t}$ to be zero  for $t\in[0,1]$. 
Notice that for $\hat{u}_{t}=tv_{1-t}$, $t\in[0,1]$, we have
$$
d\hat{u}_{t}=\big[D_{i}(\hat{a}^{ij}_{t}D_{j}\hat{u}_{t}+
\hat{f}^{i}_{t})+\hat{f}^{0}_{t}-\lambda \hat{u}_{t}\big]\,dt.
$$
Moreover, $\hat{u}_{0}=0$, $\hat{u}_{1}=u_{0}$, and $\hat{u}_{t}$
is $\hat{\cF}_{t}$-adapted. Therefore, naturally we define
$\hat{u}_{t}=u_{t-1}$ for $t\geq 1$.

It is easy to see
that if we construct the operators  
$\hat{L}_{t}$ and $\hat{\Lambda}^{k}_{t}$ from the 
coefficients with hats, then
$$
d\hat u_{t}=(\hat L_{t}\hat u_{t}-\lambda 
\hat u_{t}+D_{i}\hat f^{i}_{t}+\hat f^{0}_{t})\,dt
+(\hat \Lambda^{k}_{t}\hat u_{t}+\hat g^{k}_{t})\,d\hat w^{k}_{t},
\quad t\leq\hat \tau.
$$
By the first part of the proof
$$
\lambda\|u\|^{2}_{\bL _{p}(\tau)}+\|Du\|^{2}_{\bL _{p}(\tau)}
\leq\lambda\|\hat u\|^{2}_{\bL _{p}(\hat\tau)}+
\|D\hat u\|^{2}_{\bL _{p}(\hat\tau)}
$$
$$
\leq N\big(\sum_{i=1}^{d}\|\hat f^{i}\|^{2}_{\bL _{p}(\hat\tau)}
+\|\hat g\|^{2}_{\bL _{p}(\hat\tau)} \big)
+N\lambda^{-1}\|\hat f^{0}\|^{2}_{\bL _{p}(\hat \tau)}
$$
$$
\leq N\big(\sum_{i=1}^{d}\|  f^{i}\|^{2}_{\bL _{p}( \tau)}
+\|  g\|^{2}_{\bL _{p}(\ \tau)} \big)
+N\lambda^{-1}\|  f^{0}\|^{2}_{\bL _{p}(\tau)}
$$
$$
+N(\|v\|^{2}_{\bL_{p}}+\|Dv\|^{2}_{\bL_{p}}).
$$
It only remains to notice that the last term is dominated
by $N\|u_{0}\|_{\tr \cW^{1}_{p}}^{2}$. The theorem is proved.

\mysection{Proof of Theorem \protect\ref{theorem 3.16.1}}
                                          \label{section 6.9.5}

Throughout this section we suppose that the
assumptions of Theorem \ref{theorem 3.16.1} are
satisfied.

Owing to Theorem \ref{theorem 3.11.1}, implying that 
the solution in $\bW^{1}_{p}(\tau)$ is unique,
and having in mind setting all data equal to zero for
$t>\tau$, we see that without loss of generality we may assume
that $\tau=\infty$. Set   
$$
\bL_{p}=\bL_{p}( \infty),\quad
\bW^{1}_{p}=\bW^{1}_{p}( \infty).
$$

We need two  auxiliary results.
 \begin{lemma}
                                         \label{lemma 9.2.1}
For any   $T,R\in(0,\infty)$ (and $\omega$),  we have  
\begin{equation}
                                                     \label{9.2.2}
\int_{0}^{T}\int_{B_{R}}(|\gb_{s}(x)|^{p'}
+|b_{s}(x)|^{p'} +c_{s}^{p'}(x)) \,dxds<\infty.
\end{equation}
\end{lemma}
This lemma is proved in the same way
as Lemma 6.1 of \cite{Kr10_2} on the basis of 
Assumptions \ref{assumption 2.7.2} (iii)
and \ref{assumption 3.11.1} and the fact that
$q\geq p'$.

The solution of our   equation will be obtained
as the weak limit of the solutions of equations
with cut-off coefficients. Therefore, the following result
is relevant.

\begin{lemma}
                                            \label{lemma 3.23.2}
Let $\phi\in C^{\infty}_{0}$,   
  $u^{m}$, $u\in \bW^{1}_{p}$, $m=1,2,...$, be such that
$u^{m}\to u$ weakly in $\bW^{1}_{p}$. 
For $m=1,2,...$   define $\chi_{m}(t)=(-m)\vee t\wedge m$,
$\gb^{i}_{mt}=\chi_{m}(\gb^{i}_{t})$,
$b^{i}_{mt}=\chi_{m}(b^{i}_{t})$, and $c_{mt}=\chi_{m}(c_{t})$.
Then  there is a sequence of bounded stopping times $\tau_{n}
\to\infty$ such that, for any $n$, the functions
\begin{equation}
                                                 \label{4.19.6}
\int_{0}^{t} (b^{i}_{ms}D_{i}u^{m}_{s},\phi)\,ds,\quad
\int_{0}^{t} (\gb^{i}_{ms}u^{m}_{s},D_{i}\phi) \,ds,\quad
\int_{0}^{t} (c_{ms}u^{m}_{s}, \phi) \,ds
\end{equation}
converge weakly in the space $\cL_{p}(\opar0,\tau_{n}\cbrk)$
as $m\to\infty$ to
\begin{equation}
                                                     \label{8.12.1}
\int_{0}^{t} (b^{i}_{s}D_{i}u_{s},\phi)\,ds,
\quad
\int_{0}^{t} (\gb^{i}_{s}u_{s},D_{i}\phi)\,ds,\quad
\int_{0}^{t} (c_{s} u_{s},\phi)\,ds,
\end{equation}
respectively.
\end{lemma}

Proof. Let $R$ be such that $\phi(x)=0$ for $|x|\geq R$.
We   take  $\tau_{n}\to\infty$ such that each of them
is bounded,
they are smaller than the ones from Lemma \ref{lemma 6.27.1},
and are such that the left hand side of \eqref{9.2.2}
with $T=\tau_{n}$ is less than $n$.

 By Corollary \ref{corollary 3.23.1} and by the
fact that (strongly) continuous operators are weakly
continuous we obtain that
$$
\int_{0}^{t} (b^{i}_{s}D_{i}u^{m}_{s},\phi) \,ds
\to
\int_{0}^{t} (b^{i}_{s}D_{i}u_{s},\phi) \,ds
$$
as $m\to\infty$ 
weakly in the space $\cL_{p}(\opar0,\tau_{n} \cbrk)$
for any $n$.
 Therefore, in what concerns the first
function  in \eqref{4.19.6}, it suffices to show that
$$
\int_{0}^{t} (D_{i}u^{m}_{s},(b^{i}_{s}-b^{i}_{ms})\phi)
 \,ds
\to0
$$
weakly in $\cL_{p}(\opar0,\tau_{n} \cbrk)$. In other words, 
it suffices to show
that for any $\xi\in \cL_{p'}(\opar0,\tau_{n} \cbrk)$
$$
E\int_{0}^{\tau_{n}}\xi_{t}
\big(\int_{0}^{t} (D_{i}u^{m}_{s},(b^{i}_{s}-b^{i}_{ms})\phi)
 \,ds
\big)\,dt\to0.
$$
This relation is rewritten as
\begin{equation}
                                                 \label{4.21.1}
E \int_{0}^{\tau_{n}}
 (D_{i}u^{m}_{s},\eta_{s}(b^{i}_{s}-b^{i}_{ms})\phi)
 \,ds\to0,
\end{equation}
where  
$$
\eta_{s}:=\int_{s}^{\tau_{n} }\xi_{t}\,dt.
$$
Observe that by the choice of $\tau_{n}$ we have
$$
E\int_{0}^{\tau_{n}}|\eta_{s}|^{p'} \int_{|x|\leq R}
|b_{s}(x)|^{p'}\,dxds
\leq E\sup_{t\leq\tau_{n}}|\eta_{s}|^{p'}
\int_{0}^{\tau_{n}} \int_{|x|\leq R}
|b_{s}(x)|^{p'}\,dxds
$$
$$
\leq nE\big(\int_{0}^{\tau_{n}}|\xi_{s}|\,ds\big)^{p'}<\infty.
$$
It follows by  the 
dominated convergence that
$ 
\eta_{s}(b^{i}_{s}-b^{i}_{ms})\phi\to0 
$ 
as $m\to\infty$
strongly in $\bL_{p'}( \tau_{n}  )$. 
By assumption   $Du^{m}\to Du$
weakly in $\bL_{p}( \tau_{n} )$. This implies
\eqref{4.21.1}. Similarly,
one proves our assertion about
the remaining functions in \eqref{4.19.6}.  The lemma is proved.

{\bf Proof of Theorem \ref{theorem 3.16.1}}.
Recall that we may assume that $ \tau=\infty$.
Since the case $p=2$ is dealt with in \cite{Kr10_1}
(under much milder assumptions),
we also assume that $p>2$.
 Define
$\gb_{mt}$, $b_{mt}$, and $c_{mt}$   as in  
 Lemma \ref{lemma 3.23.2} and consider equation
\eqref{2.6.4} with $\gb_{mt}$, $b_{mt}$, and $c_{mt}$  
in place of $\gb_{t}$, $b_{t}$, and $c_{t}$, respectively.
Obviously, $\gb_{mt}$, $b_{mt}$, and $c_{mt}$  satisfy Assumption
\ref{assumption 3.11.1} with the same $\gamma_{b}$ and $K$
as $\gb_{t}$, $b_{t}$, and $c_{t}$ do. By Theorem \ref{theorem 3.11.1}
and the method of continuity for  
$\lambda\geq\lambda_{0}(d,\delta,p,\kappa,
\rho_{0},\rho_{1},K)$ there exists
a unique solution $u^{m} $ of the modified
equation on $\bR$.

By Theorem \ref{theorem 3.11.1}  we also have
$$
\|u^{m}\|_{\bL_{p} }+\|Du^{m}\|_{\bL_{p} }\leq N,
$$
where $N$ is independent of $m$. Hence the sequence of functions
$u^{m} $ is bounded in the   space 
$\bW^{1}_{p}$ and consequently has a weak limit 
point $u\in \bW^{1}_{p}$. For simplicity of presentation
we assume that the whole sequence $u^{m}  $
converges weakly to $u$. 

Take a $\phi\in C^{\infty}_{0}$. Then
by Lemma \ref{lemma 3.23.2} for appropriate $\tau_{n}$ we have
that the functions \eqref{4.19.6}  converge to \eqref{8.12.1} weakly in
$\cL_{p}(\opar0,\tau_{n}\cbrk)$ as $m\to\infty$
for any $n$.  
Owing to \eqref{9.3.3} and the fact that bounded linear
operators are weakly continuous, the stochastic
terms in the equations for $u^{m}_{t}$
also   converge   weakly in
$\cL_{p}(\opar0,\tau_{n}\cbrk)$ as $m\to\infty$
for any $n$.  
Obviously,
the same is true for $(u^{m}_{ t},\phi)\to(u_{t},\phi)$
and the remaining terms entering 
the equation for 
$ u^{m}_{t}$.
Hence, by passing to the weak limit in the equation
for $u^{m}_{ t}$ we see that for any
$\phi\in C^{\infty}_{0}$ equation
\eqref{3.16.7} holds  for {\em almost any\/} $(\omega,t)$. 

Until this moment Assumption \ref{assumption 3.16.1}
was not needed. We will need it in order to be able
to apply Theorem 3.1 of \cite{KR} and find
an appropriate modification of $u_{t}$.

Take 
a $\psi\in C^{\infty}_{0}$ and observe that
$u\psi\in\bW^{1}_{2}(T)$ and $g\psi\in\bL_{2}(T)$
for any $T\in(0,\infty)$
which implies that
$$
m^{\psi}_{t}:=u_{0}\psi+
\sum_{k=1}^{\infty}\int_{0}^{t} \psi
(\Lambda^{k}_{s}u_{s}
+g^{k}_{s} )\,dw^{k}_{s}
$$
is well defined as an $\cL_{2}$-valued continuous
martingale such that for any $\phi\in \cL_{2}$
with probability one
\begin{equation}
                                              \label{1.26.1}
(m^{\psi}_{t},\phi)=(u_{0}\psi,\phi)+
\sum_{k=1}^{\infty}\int_{0}^{t}\big( \psi
(\Lambda^{k}_{s}u_{s}
+g^{k}_{s} ),\phi\big)\,dw^{k}_{s}
\end{equation}
for all $t\in\bR_{+}$.

 Notice that for any $\phi\in C^{\infty}_{0}$
\begin{equation}
                                                 \label{9.3.4}
(u_{t}\psi,\phi)=\int_{0}^{t}(u^{*}_{s}\psi ,\phi)\,ds+
(m^{\psi}_{t},\phi)
\end{equation}
for almost all $(\omega,t)$, where $u^{*}_{s}$
is a function with values in the space of distributions
on $\bR^{d}$ defined by
$$
u^{*}_{s}= L_{s}u_{s}-\lambda u_{s}+D_{i}f^{i}_{s}+f^{0}_{s} 
$$
(see Remark \ref{remark 9.1.1}).

Next, take an $R\in(0,\infty)$ and let $W^{-1}_{p'}(B_{R})$
 denote the dual space for 
$$
\WO^{1}_{p}(B_{R})
:= W^{1}_{p}(B_{R})\cap\{v:v|_{\partial B_{R}}=0\}.
$$
Estimate \eqref{1.15.4} combined with the facts that, $p'<p$
and that one can cover $B_{R}$ with finitely many balls
of radius $\rho_{1}$ shows that for any
$\phi\in C^{\infty}_{0}(B_{R})$
$$
|(D_{i}(\gb^{i}_{s}u_{s}),\phi)|\leq N\big(1+
\int_{B_{R+1}}|\gb_{s}|\,dx\big)
\|u_{s}\|_{W^{1}_{p}} \|\phi\|_{W^{1}_{p}} ,
$$
where $N$ is independent of $\omega,s, u_{s},\phi$.
Due to the arbitrariness of $\phi$ 
and the fact that
$C^{\infty}_{0}(B_{R})$ is dense in $\WO^{1}_{p}(B_{R})$
we conclude that
(for almost all $(\omega,s)$) we have
$D_{i}(\gb^{i}_{s}u_{s})\in W^{-1}_{p'}(B_{R})$ and
$$
\|D_{i}(\gb^{i}_{s}u_{s}) \|_{W^{-1}_{p'}(B_{R})}
\leq
N\big(1+
\int_{B_{R+1}}|\gb_{s}|\,dx\big)
\|u_{s}\|_{W^{1}_{p}}.
$$
Here the right-hand side is locally summable on $\bR_{+}$
 to the power  $p'$ (a.s.) owing to Assumption \ref{assumption 3.16.1},
H\"older's inequality, and the fact that $u\in \bW^{1}_{p}$.
Similar statements are true for $b^{i}_{s}D_{i}u_{s}$, $c_{s}u_{s}$,
and $u^{*}_{s}$.

Now, since $u\psi\in\cL_{p}(\bR_{+},\WO^{1}_{p}(B_{R}))$
and $\WO^{1}_{p}(B_{R})$ is dense in $\cL_{2}(B_{R})$,
 by Theorem 3.1 of \cite{KR} we get that there exist an event
$\Omega^{\psi}$ of full probability and a continuous
$\cL_{2}(B_{R})$-valued  $\cF_{t}$-adapted process $u^{\psi}_{t}$
such that $u^{\psi}_{t}=u_{t}\psi$
as $\cL_{2}(B_{R})$-valued functions for almost all $(\omega,t)$
and for any  $\omega\in\Omega^{\psi}$, $t\in\bR_{+}$, and
$\phi\in C^{\infty}_{0}(B_{R})$ we have 
\begin{equation}
                                                 \label{9.3.7}
(u^{\psi}_{t},\phi)=\int_{0}^{t}(u^{*}_{s}\psi,\phi)\,ds+
(m^{\psi}_{t},\phi).
\end{equation}

Take a $\psi\in C^{\infty}_{0}$ such that $\psi(x)=1$ for
$|x|\leq1$ and for $k=1,2,...$ define
$\psi_{k}(x)=\psi(x/k)$ and 
$$
\Omega'=\bigcap_{k=1}^{\infty}\Omega^{\psi_{k}}.
$$
Clearly, $P(\Omega')=1$. We will further reduce
$\Omega'$ in the following way. 
Obviously (see \eqref{1.26.1}), if $\psi',\psi''\in C^{\infty}_{0}$
and $\psi'=\psi''$ on $B_{R}$ and $\phi\in\cL_{2}$
is such that $\phi=0$ outside $B_{R}$, then
with probability one we have 
$(m_{t}^{\psi'},\phi)=
(m_{t}^{\psi''},\phi)$  for all $t$.

Let $\Phi$ be the union over $n=1,2,...$ of countable
subsets of $C^{\infty}_{0}(B_{n})$ each of which 
everywhere dense in  $\cL_{2}(B_{n})$.
For $\phi\in C^{\infty}_{0}$ denote $d(\phi)$
the smallest radius of the balls centered at the origin
containing the support of $\phi$. Then by the above
for $\phi\in C^{\infty}_{0}$ the events
$$
\Omega(\phi)=\{\omega\in\Omega:
(m_{t}^{\psi_{k}},\phi)=
(m_{t}^{\psi_{j}},\phi),\quad\forall t\in\bR_{+},
 k,j\geq d(\phi) \},
$$
$$
\Omega''=\Omega'\bigcap\bigcap_{\phi\in\Phi}\Omega(\phi)
$$
have probability one. Since $m^{\psi}_{t}$ are $\cL_{2}$-valued
and $\Phi\cap C^{\infty}_{0}(B_{n})$ is dense in  
$\cL_{2}(B_{n})$, we have that for
$\omega\in\Omega''$, $t\in\bR_{+}$, and any $\phi\in
C^{\infty}_{0}(B_{n})$ it holds that
$$
(m_{t}^{\psi_{k}},\phi)=
(m_{t}^{\psi_{j}},\phi)
$$
as long as $i,j\geq n$.

Then \eqref{9.3.7} implies that for any 
$\omega\in \Omega''$, $t\in\bR{_+}$, and
$\phi\in C^{\infty}_{0}(B_{n})$ 
 we have $(u_{t}^{\psi_{j}},\phi)
=(u_{t}^{\psi_{k}},\phi)$ for all  
$j,k\geq n$. In particular, 
for any 
$\omega\in \Omega''$, $t \in\bR_{+}$, $n=1,2...$ it holds that
$u_{t}^{\psi_{j}}=u_{t}^{\psi_{k}}$
as distributions on $B_{n}$ for  
$j,k\geq n$ and there exists a distribution $\bar{u}_{t}$
on $\bR^{d}$ such that $\bar{u}_{t}=u_{t}^{\psi_{k}}$
on $B_{n}$ for all   $k\ge n$. Since
$u_{t}^{\psi_{k}}=u_{t}\psi_{k}$ for almost all
$(\omega,t)$, we have that $\bar{u}_{t}=u_{t}$
(as distributions on $\bR^{d}$) for almost all $(\omega,t)$.
The inclusion $u\in\bW^{1}_{p}$ now yields
$\bar{u}\in\bW^{1}_{p}$. 

It also 
follows from \eqref{9.3.7} that if $\omega\in\Omega''$,
$t\in\bR_{+}$, and $\phi\in C^{\infty}_{0}$ is such that
$\phi=0$ outside $B_{n}$, then for any $j\geq n$
$$
(\bar{u}_{t},\phi)=(u_{t}^{\psi_{j}},\phi)
=\int_{0}^{t}(L_{s}u_{s}-\lambda u_{s}+D_{i}f^{i}_{s}
+f^{0}_{s},\phi)\,ds+(m^{\psi_{j}}_{t},\phi).
$$
By having in mind \eqref{1.26.1} we conclude that 
for any $\phi\in C^{\infty}_{0}$ with probability
one for all $t\in\bR_{+}$
$$
(\bar{u}_{t},\phi)
=(u_{0},\phi)+
\int_{0}^{t}(L_{s}u_{s}-\lambda u_{s}+D_{i}f^{i}_{s}
+f^{0}_{s},\phi)\,ds
$$
$$
+\sum_{k=1}^{\infty}\int_{0}^{t}(\Lambda^{k}_{s}
u_{s}+g^{k}_{s},\phi)\,dw^{k}_{s}.
$$
Now it only remains to observe that
since $\bar{u}_{s}=u_{s}$ for almost all $(\omega,s)$,
we can replace  $u_{s}$ with $\bar{u}_{s}$ in the above equation.
The theorem is proved.

\mysection{It\^o's formula for the product
of two processes of class $\cW^{1}_{2,loc}(\tau)$}

                                             \label{section 1.1.2}

The results of this section will be used in a few places below,
in particular,
in the proof of Lemma \ref{lemma 10.8.2}.
Recall that the spaces $\cW^{1}_{p}(\tau)$ are
introduced in Definition \ref{definition 3.16.1}.
\begin{theorem}
                                     \label{theorem 10.14.1}
Let $\tau$ be a stopping time and let
$u,\tilde{u}$, $f^{j}$, $\tilde{f}^{j}$,
$g=(g^{1},g^{2},...)$,
$\tilde{g}=(\tilde{g}^{1},\tilde{g}^{2},...)$
be some functions such that for any $\phi\in
C^{\infty}_{0}$ we have
$\phi u,\phi\tilde{u}\in\cW^{1}_{2}(\tau)$,
$\phi f^{j},\phi\tilde{f}^{j}\in \bL_{2}(\tau)$,
$j=0,...,d$, and $\phi g ,
\phi\tilde{g} \in \bL_{2}(\tau)$. Assume that
in the sense of generalized functions
$$
du_{t}=(D_{i}f^{i}_{t}+f^{0}_{t})\,dt+g^{k}_{t}\,dw^{k}_{t},
\quad
d\tilde{u}_{t}=(D_{i}\tilde{f}^{i}_{t}+
\tilde{f}^{0}_{t})\,dt+\tilde{g}^{k}_{t}\,dw^{k}_{t},
\quad t\leq\tau.
$$

Then
$$
d(u_{t}\tilde{u}_{t})=\big[
\tilde{u}_{t}(D_{i}f^{i}_{t}+f^{0}_{t})
+u_{t}(D_{i}\tilde f^{i}_{t}+\tilde f^{0}_{t})+h_{t}\big]\,dt
$$
$$
+
(\tilde{u}_{t}g^{k}_{t}+u_{t}\tilde{g}^{k}_{t})\,dw^{k}_{t},
\quad t\leq\tau,
$$
where $h_{t}:=
(g _{t},\tilde{g} _{t}
)_{\ell_{2}}$,
in the sense of generalized functions, that is,
for any $\phi\in C^{\infty}_{0}$, with probability one,
$$
 (u_{t\wedge\tau}\tilde{u}_{t\wedge\tau},\phi)=
(u_{0}\tilde{u}_{0},\phi)
+\int_{0}^{t}I_{s\leq\tau}
  (\tilde{u}_{s}g^{k}_{s}+u_{s}\tilde{g}^{k}_{s}
,\phi)\,dw^{k}_{s}
$$
\begin{equation}
                                          \label{10.16.1}
+\int_{0}^{t}I_{s\leq\tau}
\big[(\tilde{u}_{s}f ^{0}_{s},\phi)
 -\big( f^{i}_{s},D_{i}(\tilde{u}_{s}\phi)\big)+
 (u_{s}\tilde{f}^{0}_{s},\phi)
 -\big( \tilde{f}^{i}_{s},D_{i}(u _{s}\phi)\big)
+(h_{s},\phi)\big]\,ds
\end{equation}
for all $t$.
\end{theorem}

Proof. To prove \eqref{10.16.1},
we only need to consider the case that $\tilde{u}=u$.
Indeed, then by writing down the stochastic differential
of $|u_{t}+\lambda\tilde{u}_{t}|^{2}$, where $\lambda$
is an arbitrary constant, and comparing the coefficients
of $\lambda$, we would come to \eqref{10.16.1}.
In other words, to prove \eqref{10.16.1},
 we need only prove that for any $\phi
\in C^{\infty}_{0}$ with probability one
$$
 (u^{2}_{t\wedge\tau} ,\phi)=
(u^{2}_{0} ,\phi)
+2\int_{0}^{t}I_{s\leq\tau}
  (u_{s}g^{k}_{s} 
,\phi)\,dw^{k}_{s}
$$
\begin{equation}
                                       \label{10.16.3}
+\int_{0}^{t}I_{s\leq\tau}
\big[2(u_{s}f ^{0}_{s},\phi)
 -2\big( f^{i}_{s},D_{i}(u_{s}\phi)\big) 
  +(|g_{s}|^{2}_{\ell_{2}},\phi)\big]\,ds
\end{equation}
for all $t$.

Next, observe that for any $\psi,\phi\in C^{\infty}_{0}$,
with probability one
$$
 (\psi u_{t\wedge\tau},\phi)= (\psi u_{0},\phi)+
\int_{0}^{t}I_{t\leq\tau}
(\psi g^{k}_{s},\phi)\,dw^{k}_{s}
$$ 
$$
+\int_{0}^{t}I_{t\leq\tau}\big[
(\psi f^{0}_{s}-f^{i}_{s}D_{i}\psi,\phi)
-(\psi f^{i}_{s},D_{i}\phi)
\big]\,dt.
$$
for all $t$. This means that
$$
d(\psi u_{t})=(\psi f^{0}_{t}-f^{i}_{t}D_{i}\psi
+D_{i}(\psi f^{i}_{t}))\,dt
+\psi g^{k}_{t} \,dw^{k}_{t},\quad t\leq\tau.
$$
By well-known results, in particular, by
It\^o's formula (see, for instance \cite{Kr09_3})
there is a set $\Omega'\subset\Omega$ of full probability
such that

(i) $\psi u_{t\wedge\tau}I_{\Omega'}$
is a continuous $\cL_{2}$-valued $\cF_{t}$-adapted function on
$[0,\infty)$;

(ii) for all $t\in[0,\infty)$ and $\omega\in\Omega'$,
It\^o's formula holds:
$$
\int_{\bR^{d}}|\psi u_{t\wedge\tau}|^{2}\,dx
=\int_{\bR^{d}}|\psi u_{0}|^{2}\,dx +2
\int_{0}^{t  }I_{s\leq\tau}
\int_{\bR^{d}} \psi^{2} u _{s}
g^{k }_{s}\,dx\,dw^{k}_{s}
$$
\begin{equation}
                                            \label{4.19.5}
+
\int_{0}^{t }I_{s\leq\tau}
\big(
\int_{\bR^{d}}\big[2 
u_{ s}f^{0}_{ s}\psi^{2}
-2 f^{i}_{ s}
D_{i}(\psi^{2} u_{ s})
+  \psi^{2}|
g_{ s }|_{\ell_{2}}^{2}
\big]\,dx\big)\,d s. 
\end{equation}

This proves \eqref{10.16.3} if we replace there $\phi$
with $\psi^{2}$. However, for any $\phi\in C^{\infty}_{0}$
one can find $\psi_{1},\psi_{2}\in C^{\infty}_{0}$
such that $\phi=\psi_{1}^{2}-\psi_{2}^{2}$. Indeed,
  one can take
 sufficiently large $N,R>0$  
and take $\psi_{1}(x)=\exp(-(R^{2}-|x|^{2})^{-1})$
for $|x|<R$ and $\psi_{1}(x)=0$ for $|x|\geq R$ and define
$\psi_{2}=(\psi_{1}^{2}-\phi)^{1/2}$.
This implies that \eqref{10.16.3} holds
for any $\phi\in C^{\infty}_{0}$ with probability one
for all $t$  and proves the theorem.
\begin{corollary}
                               \label{corollary 09.10.29.1}

Let $u,f,g$ be as in Theorem \ref{theorem 10.14.1},
let a nonrandom $\psi\in W^{1}_{2}$, and let
a random process $x_{t}$ be given as
$$
x_{t}=\int_{0}^{t}\sigma^{k}_{s}\,dw^{k}_{s}
+\int_{0}^{t}b_{s}\,ds
$$
for some predictable $\bR^{d}$-valued functions
 $\sigma^{k}_{t}$ and $b_{t}$ such that
$$
E\int_{0}^{\tau}\big(\sum_{k}|\sigma^{k}_{t}|^{2}
+|b_{t}|\big)\,dt<\infty.
$$
Then in the sense of generalized functions
$$
d(u_{t}\psi_{t})=\big[D_{i}(u_{t}a^{ij}_{t}D_{j}\psi_{t})
-a^{ij}_{t}(D_{i}u_{t})D_{j}\psi_{t}
+u_{t}b_{t}^{i}D_{i}\psi_{t}+D_{i}(\psi_{t}f^{i}_{t})
$$
$$
-f^{i}_{t}D_{i}\psi_{t}+\psi_{t}f^{0}_{t}
+g^{k}_{t}\sigma^{ik}_{t}D_{i}\psi_{t}\big]\,dt
+\big[\psi_{t}g^{k}_{t}+u_{t}\sigma^{ik}_{t}D_{i}
\psi_{t}\big]\,dw^{k}_{t},\quad t\leq\tau,
$$
where 
$\psi_{t}(x)=\psi(x+x_{t})$ and $2a^{ij}_{t}
=\sigma^{ik}_{t}\sigma^{jk}_{t}$.
\end{corollary}

Indeed,  observe that by It\^o's formula
and the stochastic Fubini theorem, for any $\phi
\in C^{\infty}_{0}$,
$$
\int_{\bR^{d}}\psi_{t\wedge\tau}\phi\,dx=
\int_{\bR^{d}}\psi(x) \phi(x-x_{t\wedge\tau})\,dx=
\int_{\bR^{d}}\psi\phi\,dx
$$
 $$
+\int_{0}^{t}I_{s\leq\tau}\int_{\bR^{d}}
\psi_{s}[a^{ij}_{s}D_{ij}\phi
-b^{i}_{s}D_{i}\phi]\,dx\,ds+
\int_{0}^{t}I_{s\leq\tau}\int_{\bR^{d}}\sigma^{ik}_{s}
\phi D_{i}\psi_{s}\,dx\,dw^{k}_{s},
$$
where the coefficient of $ds$ equals
$$
\int_{\bR^{d}}
\phi[a^{ij}_{s}
 D_{ij}\psi_{s}+
 b^{i}_{s} D_{i}\psi_{s}]\,dx.
$$
Furthermore, for instance,
$$
E\int_{0}^{\tau}\int_{\bR^{d}}\sum_{k} |\sigma^{ik}_{s}
  D_{i}\psi_{s}|^{2}\,dx\,ds
\leq E\int_{0}^{\tau}\int_{\bR^{d}}\sum_{k} |\sigma^{ k}_{s}
 |^{2}| D\psi_{s}|^{2}\,dx\,ds
$$
$$
=\int_{\bR^{d}}| D\psi |^{2}\,dx\,
E\int_{0}^{\tau}\sum_{k} |\sigma^{ k}_{s}
 |^{2}\,ds<\infty.
$$
It follows that $\psi_{\cdot}\in\cW^{1}_{2}(\tau)$ and
$$
d\psi_{t}=\big[D_{i}(a^{ij}_{t}D_{j}\psi_{t})
+b^{i}_{t}D_{i}\psi_{t}\big]\,dt
+\sigma^{ik}_{t}D_{i}\psi_{t}\,dw^{k}_{t}
$$
in the sense of generalized functions, so that the desired 
result follows from Theorem \ref{theorem 10.14.1}.

\mysection{Kalman-Bucy filter}
                                \label{section 11.8.1}

We take a $T\in(0,\infty)$ and on $[0,T]$ consider a
$d_{1}$-dimensional two component process
 $z_{t}=(x_{t},y_{t})$ with $x_{t}$ being $d$-dimensional
and $y_{t}$ 
$(d_{1}-d)$-dimensional. We assume that $z_{t}$ is a
diffusion process
 defined as a solution of the system
\begin{equation}\begin{split}         
                                           \label{eq3.2.14}
  & dx_{t}=b(t,z_{t}) dt+\theta
(t,y_{t})dw_{t}, \\  & dy_{t}=B(t,z_{t}) dt+\Theta
(t,y_{t})dw_{t}
\end{split}\end{equation}
with some initial data.

\begin{assumption}
                                           \label{asm3.2.15}
The functions $b$, $\theta $, $B$ and $\Theta $ are
 Borel measurable  
functions of $(t,z)$ and $(t,y)$
as appropriate and $\theta$ and $\Theta $
are bounded and satisfy the Lipschitz condition
with respect to $y$ with a constant independent
of $t$. 
We have
$$
b(t,z)=x^{*}\dot{b}(t,y)+b(t,0,y),\quad
B(t,z)=x^{*}\dot{B}(t,y)+B(t,0,y),
$$
where $\dot{b}$ and $\dot{B}$ are bounded matrix-valued functions
of appropriate dimensions, $b(t,0)$ and $B(t,0)$
are bounded, and $\dot{b}(t,y)$, $\dot{B}(t,y)$,
$b(t,0,y)$, and $B(t,0,y)$ satisfy the Lipschitz condition
with respect to $y$ with a constant independent
of $t$. 
 
\end{assumption}

In the rest of the article we use the notation
$$
D_{i}=\frac{\partial}{\partial x^{i}},\quad D_{ij}=D_{i}D_{j}
$$
only for $i,j=1,...,d$.

\begin{remark}
Note that
\begin{equation}
                                            \label{9.19.1}
 \dot{b}^{ij}(t,y)
=D_{i}b^{j} (t,z),\quad\dot{B}^{ij}(t,y)
=D_{i}B^{j} (t,z).
\end{equation}
\end{remark}

Set
\begin{equation}
                                      \label{eq3.2.19.4}
\check{\theta } (t,y) =\begin{pmatrix}
\theta(t,y)\\
\Theta(t,y)
\end{pmatrix},\quad
\check{a} (t,y) =\frac{1}{2}
\check{\theta }\check{\theta }^{*}(t,y),\quad
\check{b} (t,z) =\begin{pmatrix}
b(t,z) \\
B(t,z)
\end{pmatrix},
\end{equation}
\begin{equation}
                                             \label{eq3.2.19.1}
\check{L} (t,z)  =   \check{a}^{ij} (t,y)
\frac{\partial^{2} }{\partial z^{i} \partial z^{j}}+
 \check{b}^{i} (t,z)\frac{\partial }{\partial z^{i}},
\end{equation}
where $t\in[0,T]$,   $z=(x,y)\in \bR^{d_{1}}$,  
and we use the summation convention over all ``reasonable''
values of repeated indices,
so that the summation  in
(\ref{eq3.2.19.1}) is performed
  for $i,j=1,...,d_{1}$.

Observe that
\begin{equation}
                                                \label{9.29.4}
dz_{t}=\check{\theta }(t,z_{t})\,dw_{t}+
\check{b}(t,z_{t})\,dt.
\end{equation}

\begin{assumption}
                                           \label{asm3.2.16}
The symmetric matrix $\check{a} (t,y)$ is uniformly
nondegenerate. In particular,
the matrix
$\Theta \Theta^{*}$ is 
invertible and 
$$
\Psi :=(\Theta
\Theta^{*} )^{-\frac{1}{2}}
$$ 
is a bounded function of $(t,y)$.
\end{assumption}

\begin{remark}
                                          \label{remark 9.19.1}

It is well known (see, for instance, \cite{Kr_10})
that in light of Assumption \ref{asm3.2.16}
the matrix 
$$
\hat{a}(t,y)=a(t,y)-\alpha(t,y)
$$
is uniformly (with respect to $(t,y)$) nondegenerate, where
$$
a =\frac{1}{2}\theta \theta^{*} ,
\quad \alpha=\tfrac{1}{2}\sigma\sigma^{*},\quad
\sigma =\theta\Theta^{*}\Psi ,
$$

\end{remark}
\begin{remark}
                                        \label{remark 1.16.1}
Everywhere below we use the stipulation 
that if we are given a function $\xi(t,x,y)$, then
we denote 
\begin{equation}
                                                 \label{1.16.1}
\xi_{t}=\xi_{t}(x)=\xi(t,x,y_{t})
\end{equation}
unless it is explicitly specified otherwise.
For instance, $
\Psi_{t}=\Psi(t,y_{t})$, $\Theta_{t}=\Theta(t,y_{t})$,
$
\sigma_{t}  =\theta_{t} \Theta^{*}_{t}\Psi_{t}$.
\end{remark}
Next we introduce a few more notation.
Let  (note the size and shape  of $\sfB$)
$$
\sfB=\Psi B,\quad
\sfB_{t}(x) =\Psi_{t}B_{t}(x)=\Psi(t,y_{t})
B(t,x,y_{t}) 
$$
 and set
\begin{equation}                             \label{eq3.2.19.2}
L_{t}( x) = a^{ij}_{t} D_{i}D_{j}+
 b^{i}_{t}(x)D_{i}\,,
\end{equation}
$$
L^{*}_{t}( x)u_{t}(x) = 
D_{i}D_{j}(
a^{ij}_{t} u_{t}(x) )
- D_{i}(b^{i}_{t}( x)u_{t}(x) )
$$
\begin{equation}                     \label{e3.2.19.2}
=D_{j}\big(
a^{ij}_{t} D_{i}u_{t}(x)
-b^{j}_{t}(x)u_{t}(x) \big),
\end{equation}
\begin{equation}                     \label{1.27.8}
\Lambda^{k }_{t}( x)u_{t}(x)  = 
 \sigma^{ik}_{t}  D_{i}u_{t}(x)+\sfB^{k}_{t}( x)u_{t}(x),
\end{equation}
\begin{equation}\label{eq3.2.19.3}
\Lambda^{k*}_{t}( x)u_{t}(x) =-\sigma^{ik}_{t} D_{i}u_{t}(x)
+ \sfB^{k}_{t}(x) u_{t}(x),
\end{equation}
where $t\in[0,T]$,   $x\in \bR^{d}$, $k=1,...,d_{1}-d$,
and as above
 we use the summation convention over all ``reasonable''
values of repeated indices,
so that the summation in (\ref{eq3.2.19.2}), (\ref{e3.2.19.2}),
\eqref{1.27.8},
and (\ref{eq3.2.19.3}) is performed for $i,j=1,...,d$ (whereas in
(\ref{eq3.2.19.1}) for $i,j=1,...,d_{1}$).

Finally, by $\cF_{t}^{y}$ we denote the completion of 
$\sigma \{ y_{s}:s\leq t\}$
with respect to $P,\cF$.
\begin{assumption}
                                     \label{assumption 9.20.1}
There exists   an $\varepsilon>0$ and
a function $Q(x)=Q(\omega,x)$ which is
  $\cF^{y}_{0}$-measurable in $\omega$, quadratic in $x$, and

(i) For all $x\in\bR^{d}$ (and $\omega$)  
$$
\varepsilon^{-1}|x|^{2}\geq  x^{i}x^{j}D_{ij}Q
\geq\varepsilon |x|^{2};
$$

(ii) We have
  $\pi_{0}e^{Q}\in \tr\cW^{1}_{p}$,
where $\pi_{0}$ is   the conditional density 
of $x_{0}$ given~$y_{0}$.
\end{assumption}

Assumption \ref{assumption 9.20.1} is satisfied, for
instance, in the classical setting of the Kalman-Bucy
filter when $\pi_{0}$ is a Gaussian density.  

\begin{theorem}
                                     \label{theorem 11.8.1}
There exists a  process
$\bar{\pi } $ on $[0,T]$  such that   

(i) $\bar{\pi }_{t}$ is $\cF^{y}_{t}$-adapted
and, for any $r\in[1,p]$,
with probability one $\bar{\pi}_{t}$
is a continuous $\cL_{r}$-valued process on $[0,T]$
and $\bar{\pi}_{0}=\pi_{0}$;

(ii) There exists an increasing
 sequence of
$\cF^{y}_{t}$-stopping times $\tau_{m}\leq T$ such that
$P(\tau_{m}=T)\to1$ and $\bar{\pi}
\in\bW^{1}_{p}(\tau_{m})$ 
for any $m$;
 
(iii) In the sense of Definition \ref{definition 3.20.01}
for any $m$
\begin{equation}
                                        \label{10.20.2}
d\bar{\pi}_{t}=\Lambda^{k*}_{t}\bar{\pi}_{t}
\,d\tilde{y}^{k}_{t}+L^{*}_{t}\bar{\pi}_{t}\,dt,\quad 
t\leq \tau_{m} ,
\end{equation}
where
$$
\tilde{y}^{k}_{t}=\int_{0}^{t}
\Psi^{kr}_{s}\,dy^{r}_{s}.
$$
Furthermore, for any $m$ and $\phi\in C^{\infty}_{0}$
we have $\bar{\pi}\phi\in\cW^{1}_{p}(\tau_{m})$;

(iv) We have 
$\bar{\pi}_{t}\geq 0$ for all $t\in[0,T]$ (a.s.),  
\begin{equation} 
                                             \label{1.27.07} 
0<\int_{\bR^{d}} 
\bar{\pi }_{t}(x)\,dx=(\bar{\pi }_{t},1)<\infty 
\end{equation}
 for 
all $t\in[0,T]$  (a.s.),
and for any $t\in[0,T]$ and real-valued,
bounded or nonnegative, (Borel) measurable function $f$ 
given on $\bR^{d}$
\begin{equation} 
                                             \label{10.13.3}
E[f(x_{t})|\cF_{t}^{y}]=
\frac{(\bar{\pi }_{t},f)}
{(\bar{\pi }_{t},1)}\quad
\text{(a.s.).}
\end{equation}
\end{theorem}

\begin{remark}
Equation \eqref{10.13.3} shows (by definition) that
$$
\pi_{t}(x):=\frac{\bar{\pi }_{t}(x) }
{(\bar{\pi }_{t},1)}
$$
is a conditional density of distribution of
$x_{t}$ given $y_{s},s\leq t$. Since, generally,
$(\bar{\pi }_{t},1)\ne1$,  one calls 
$\bar{\pi }_{t} $
an unnormalized conditional density of distribution of
$x_{t}$ given $y_{s},s\leq t$. Thus, 
Theorem \ref{theorem 11.8.1} allows us to characterize
the conditional density and being combined with Theorem
\ref{theorem 1.12.1} allows us to obtain fine regularity
properties of it.
\end{remark}
The following result is obtained by repeating what is said
after Theorem~\ref{theorem 1.12.1} and taking into account that
with probability one $\tau_{m}=T$ for all large $m$.
\begin{theorem}
                              \label{theorem 1.17.1}

(ii) For any $\phi\in C^{\infty}_{0}$ the process
 $\bar{\pi }_{t }\phi$ is continuous on $[0,T]$
as an $\cL_{p}$-valued process (a.s.);

(ii) If $p>2$  and we have two numbers
$\alpha$ and $\beta$ such that
$$
\frac{2}{p}<\alpha<\beta\leq1,
$$
 then
for any $\phi\in C^{\infty}_{0}$ (a.s.)
$$
\bar{\pi }\phi\in C^{\alpha/2-1/p}([0,T], H^{1-\beta}_{p}).
$$
In particular, if $p>d+2$, then 

(a) for any $\varepsilon
\in(0,\varepsilon_{0}]$, with
$\varepsilon_{0}=1-(d+2)/p$, (a.s.) for any
$t\in[0,T]$ we have $\bar{\pi }_{t}\phi\in C^{\varepsilon_{0}-
\varepsilon}(\bR^{d})$ and the norm of $\bar{\pi }_{t}\phi$
in this space is bounded as a function of $t$;

(b) for any $\varepsilon$ as in (a) (a.s.)
for  any $x\in\bR^{d}$
we have $\bar{\pi }_{\cdot}(x)\phi(x)\in C^{(\varepsilon_{0}-
\varepsilon)/2}([0,T])$ and the norm of 
$\bar{\pi }_{\cdot}(x)\phi(x)$
in this space is bounded as a function of $x$.

\end{theorem}

In the general filtering  theory equation \eqref{10.20.2}
is known as Zakai's equation. From the point of view of
the Sobolev space theory of SPDEs the most unpleasant feature
of \eqref{10.20.2} in our particular case is the presence
of $\sfB^{k}_{t}(x)\bar{\pi}_{t}$ in the stochastic term with 
$\sfB^{k}_{t}(x)$ which is unbounded in $x$. However, in the theory
of linear PDEs it was observed that if an equation has a zeroth
order term and we know a particular nonzero solution, then the ratio
of the unknown function and this particular solution
satisfies an equation without zeroth order 
term (cf. \eqref{9.29.01}).

The way to find a particular solution of \eqref{10.20.2}
is suggested by filtering theory. Imagine that $\check{b}$
is affine with respect to $z$ and $\check{\theta}$ is independent
of $z$. Then as easy to see $z_{t}$ is a Gaussian process and hence
the conditional density of $x_{t}$ given $y_{s}$, $s\leq t$,
is Gaussian, that is, its logarithm is a quadratic function in $x$.
Therefore, we were looking for a particular solution as
$e^{-Q_{t}(x)}$, where $Q_{t}(x)$ is   a quadratic function
with respect to $x$, and 
finding the equation for $Q_{t}(x)$ (see \eqref{11.9.1})
  was pretty 
straightforward.

After we ``kill'' the zeroth-order term  our equation
falls into the scheme of Section \ref{section 2.15.1}
 even though it still
 has growing first order coefficients in the deterministic
part of the equation. Finding $\bar{\pi}_{t}$ in the described way
allows us to follow the scheme suggested in \cite{KZ}
thus avoiding using filtering theory.
However, we still encounter an additional difficulty
that certain exponential martingales may not have
moments of order $>1$, unlike the situation in \cite{KZ},
and, to prove that they are martingales indeed, we use
the Liptser-Shiryaev theorem (see \cite{LS}).
This way of proceeding was used by Liptser in \cite{Li}
(see also \cite{LS})   while treating filtering
problem for the  so-called
 conditionally Gaussian processes. 

Finding a particular solution of \eqref{10.20.2}
is based on the following lemma which is probably well known.
We give its proof in the end of Section \ref{section 1.1.3}
just for completeness. Set  
\begin{equation}
                                          \label{9.10.4}
\dot{\sfB}_{t}=\dot{B}_{t}\Psi_{t}.
\end{equation}
\begin{lemma}
                                           \label{lemma 11.1.1}

The following system of equations about $d\times d$-symmetric
matrix-valued process $W_{t}$, $\bR^{d}$-valued process
$V_{t}$, and real-valued process $U_{t}$
\begin{equation}
                                              \label{9.12.1}
\frac{d}{dt}W_{t}= (\dot{\sfB}_{t}\sigma_{t}^{*}
-\dot{b}_{t} 
)W_{t}+W^{*}_{t}(\sigma _{t}\dot{\sfB}_{t}^{*}
-\dot{b}^{*}_{t})
-2W^{*}_{t}\hat{a}_{t}W_{t}+\dot{\sfB}_{t}
 \dot{\sfB}_{t} ^{*},
\end{equation}
$$
dV_{t}=-(W_{t}\sigma_{t}+\dot{\sfB}_{t}) \,d\tilde{y}_{t}
$$
\begin{equation}
                                              \label{9.12.2} 
+[(\dot\sfB_{t}\sigma^{*}_{t}-\dot{b}_{t} )V_{t}
-2W_{t}\hat{a}_{t}V_{t}
+W_{t}(\sigma_{t}\sfB_{t}(0)-b_{t}(0) )
+\dot{\sfB}_{t}\sfB_{t}(0)
]\,dt,
\end{equation}
$$
dU_{t}=-(V^{*}_{t}\sigma_{t}+\sfB^{*}_{t}(0)) \,d\tilde{y}_{t}
$$
\begin{equation}
                                              \label{9.15.1}
+[a^{ij}_{t}W^{ij}_{t}+
V^{*}_{t}(\sigma_{t}\sfB_{t}(0)-b_{t}(0) )
-V^{*}_{t}\hat{a}_{t}V _{t}
+\tfrac{1}{2}|\sfB_{t}(0)|^{2}+\tr\dot{b}_{t}]\,dt,
\end{equation}
has a unique $\cF^{y}_{t}$-adapted
 solution with initial conditions
$W_{0}^{ij}=D_{ij}Q$, $V_{0}^{i}=D_{i}Q(0)$, 
$U_{0}=Q(0) $. Furthermore, $\varepsilon_{1}^{-1}(\delta^{ij})
\geq W_{t}\geq\varepsilon_{1}(\delta^{ij})$
on $[0,T]$, where $\varepsilon_{1}>0$ is a constant
independent of $\omega$ and $t$
(depending on $T$ among other things).
\end{lemma}

Observe that the coefficients 
in \eqref{9.15.1} are  independent of $x$.

\begin{remark}
Set
\begin{equation}
                                                 \label{11.11.1}
Q_{t}(x)=\tfrac{1}{2}W^{ij}_{t}x^{i}x^{j}+
V^{i}_{t}x^{i}+U_{t}.
\end{equation}
Then by using It\^o's formula one easily checks that
for any $x\in\bR^{d}$
$$
dQ_{t}(x)=-(\sigma^{ik}_{t}D_{i}Q_{t}(x)+
\sfB^{k}_{t}(x)
)\,d\tilde{y}^{k}_{t} +\big[ a^{ij}_{t}D_{ij}Q_{t}(x)+
D_{i}b^{i}_{t}
$$
\begin{equation}
                                                     \label{11.9.1}
+(\sigma^{ik}_{t}\sfB^{k}_{t}(x)-b^{i}_{t}(x))
 D_{i}Q_{t}(x) -
\hat{a}^{ij}_{t} 
(D_{i}Q_{t}(x))D_{j}Q_{t}(x)
 +\tfrac{1}{2}|\sfB_{t}(x)|^{2}\big]\,dt 
\end{equation}
and $\eta_{t}=e^{-Q_{t}}$ satisfies
\begin{equation}
                                             \label{9.10.5}
d\eta_{t}(x)= \Lambda^{r*}_{t}\eta_{t}(x)\,d\tilde{y}^{r}_{t}
+ L^{*}_{t}  \eta_{t}(x) \,dt.
\end{equation}

By the way,  $Q_{t}(x)$ is
a unique $\cF^{y}_{t}$-adapted function 
  depending quadratically
on $x$, satisfying \eqref{11.9.1}, and such that $Q_{0}=Q$.
Indeed, uniqueness follows from the fact that $D_{ij}Q_{t}$,
$D_{i}Q_{t}(0)$, and $Q_{t}(0)$ are easily shown to
satisfy \eqref{9.12.1}, \eqref{9.12.2}, and \eqref{9.15.1},
respectively.
\end{remark}

Our method also allows us to derive the classical
equations for the Kalman-Bucy filter. 
\begin{theorem}
                                     \label{theorem 11.11.1}
Replace requirement (ii) in Assumption \ref{assumption 9.20.1}
with the assumption that $\pi_{0}=e^{-Q}$.
 Then for any $t$ (a.s.) we have
$\pi_{t}(x)=C_{t}e^{-Q_{t}(x)}$,
 where $c_{t}$ is a normalizing process obtained
 from the condition that
$$
C_{t}\int_{\bR^{d}}e^{-Q_{t}(x)}\,dx=1.
$$

\end{theorem}

This theorem is proved in Section \ref{section 1.1.4}.

\begin{remark}
                                   \label{remark 1.17.1}

After just completing the square and finding the 
stochastic differential
 of the remaining term we find  that
$$
Q_{t}(x)=\tfrac{1}{2}|W^{1/2}_{t}x+W^{-1/2}_{t}V_{t}|^{2}
+\int_{0}^{t}(V^{*}_{s}W^{-1}_{s}\dot{\sfB}_{s}
-\sfB^{*}_{s}(0))\,d\tilde{y}_{s}
$$
\begin{equation}
                                                     \label{11.1.3}
+\tfrac{1}{2}\int_{0}^{t}|\dot{\sfB}^{*}_{s}W_{s}^{-1}
V_{s}-\sfB_{s}(0)|^{2}\,ds+A_{t},
\end{equation}
with   a bounded  
  on $\Omega\times [0,T]$ function
$$
A_{t}:=\int_{0}^{t}[a^{ij}_{s}W^{ij}_{s}+\tr \dot{b}_{s}
-\tfrac{1}{2}\|W^{1/2}_{s}\sigma_{s}
+W^{-1/2}_{s}\dot{\sfB}_{s}\|^{2}
]\,ds,
$$
where for a matrix $u$ we use the notation $\|u\|^{2}
=\tr uu^{*}$.
This shows that in the situation of Theorem \ref{theorem 11.11.1}
$$
\bar{x}_{t}:=E(x_{t}\mid\cF^{y}_{t})= \int_{\bR^{d}}x\pi_{t}(x)\,dx
=-W_{t}^{-1}V_{t},
$$
$$
\Sigma_{t}:=E\big((x_{t}-\bar{x}_{t})(x _{t}-\bar{x}_{t})^{*}
\mid\cF^{y}_{t}\big)=W_{t}^{-1}
$$
and allows one to derive the classical Kalman-Bucy
equations for $\bar{x}_{t}$ and $\Sigma_{t}$
from \eqref{9.12.1} and \eqref{9.12.2}.

\end{remark}

\mysection{An auxiliary function}

                                         \label{section 1.1.3}

The assumptions from Section \ref{section 11.8.1}
are supposed to hold. 
Set
$$
\hat{b}^{i}_{t}(x)=
\sigma^{ik}_{t}\sfB^{k}_{t}(x)-   
2\hat{a}^{ij}_{t} 
D_{j}Q_{t}(x).
$$

\begin{theorem}
                                      \label{theorem 9.19.01}
The  equation
\begin{equation}
                                              \label{9.29.01}
d\hat{\pi}_{t}=
-\sigma^{ik}_{t}D_{i}\hat{\pi}_{t}\,
d\tilde{y}^{k}_{t}
+\big[ a_{t}^{ij}D_{ij}
 \hat{\pi }_{t}-b^{i}_{t}D_{i}\hat{\pi}_{t}
 +\hat{b}^{i}_{t}D_{i}\hat{\pi}_{t} 
   \big]\,dt
,\quad t\leq T ,
\end{equation}
with initial data $\hat{\pi}_{0}=e^{Q}\pi_{0}$
has a unique solution  in the sense of 
Definition~\ref{definition 3.20.01}.

\end{theorem}  

This theorem is a direct consequence of
Remark \ref{remark 9.29.2} and Theorem \ref{theorem 3.16.1}
since the coefficients $b$ and $\hat{b}$
in \eqref{9.29.01}
are affine functions of $x$ and have {\em bounded\/}
derivatives in $x$.

\begin{lemma}
                                     \label{lemma 10.30.1}
Almost surely $\hat{\pi}_{t}$ is a continuous
$\cL_{p}$-valued process on $[0,T]$.
Furthermore, $G _{t}\|\hat{\pi}_{t}\|^{p}_{\cL_{p}}$
is a decreasing function of $t$ (a.s.), where $G_{t}$
is a bounded function on $\Omega\times[0,T]$ defined by
$$
G_{t}:=\exp\int_{0}^{t}(
 D_{i}\hat{b}^{i}_{s}-D_{i}b^{i}_{s})\,ds
=\exp\int_{0}^{t}\tr(\sigma_{s}\dot{\sfB}^{*}_{s}
-\hat{a}_{s}W_{s}-\dot{b}_{s})\,ds.
$$
In particular, on the set where $\tau:=
T\wedge\inf\{t\geq0:
\|\hat{\pi}_{t}\|_{\cL_{p}}=0\}<T$ we have
$\|\hat{\pi}_{t}\|_{\cL_{p}}=0$ for $\tau\leq t\leq T$
(a.s.).
\end{lemma}  

Proof. Set
$$
\xi^{i}_{t}=\int_{0}^{t}
\sigma^{ik}_{s}\,d\tilde{y}^{r}_{s},\quad\xi_{t}=
(\xi^{i}_{t}),\quad
\tau_{m}=T\wedge\inf\{t\geq0:|z_{t}|+|\xi_{t}|\geq m\}.
$$ 
The purpose to stop $z_{t}$ is that on $\opar0,\tau_{m}\cbrk$
we have
$$
|\sigma_{t}\sfB_{t}(x)|+|b_{t}(x)|+
|\hat{b}_{t}(x)|\leq N(1+|x|),
$$
where the constant $N$ is independent of $\omega,t,x$.
Why we also stop $\xi_{t}$ will become clear later.  

Observe that for any $\psi\in C^{\infty}_{0}$
the process $\psi\hat{\pi}_{t}$ satisfies an equation
obtained  by multiplying through \eqref{9.29.01} 
by $\psi$. Then after writing $\psi D_{i}(a^{ij}_{t}
D_{j}\hat{\pi}_{t})$ as $ D_{i}(\psi a^{ij}_{t}
D_{j}\hat{\pi}_{t})-   a^{ij}_{t}
(D_{j}\hat{\pi}_{t})D_{i}\psi$ 
and noting that the other coefficients multiplied by $\psi$
are bounded functions on $\opar0,\tau_{m}\cbrk\times\bR^{d}$
we see that
\begin{equation}
                                          \label{11.1.1}
\psi\hat{\pi}_{t}\in \cW^{1}_{p}(\tau_{m})
\end{equation}
for any $m$. It follows from \cite{Kr09_3} that
with probability one $\psi\hat{\pi}_{t\wedge\tau_{m}}$
 is a continuous $\cL_{p}$-valued process and since,
for each $\omega$,
$\tau_{m}=T$ if $m$ is sufficiently large,
with probability one $\psi\hat{\pi}_{t }$
 is a continuous $\cL_{p}$-valued process on $[0,T]$
for any $\psi\in C^{\infty}_{0}$.

Then take a nonnegative radially symmetric
and radially decreasing function
 $\phi\in C^{\infty}_{0}$ such that $|D\phi|\leq1$,
 introduce $\phi^{n} (x)
=\phi(x/n)$, $n=1,2,...$,
$$
\phi^{n}_{t}(x)=\phi^{n} (x-\xi_{t})
$$
and use Corollary \ref{corollary 09.10.29.1}
with $\tau_{m}$ in place of $\tau$
(recall \eqref{11.1.1}). Then we find
$$
d(\hat{\pi}_{t}\phi^{n}_{t})=-\sigma^{ik}_{t}D_{i}(
\hat{\pi}_{t}\phi^{n}_{t})\,d\tilde{y}^{k}_{t}
+\phi^{n}_{t} ( \hat{b}^{i}_{t}-b^{i}_{t} ) D_{i}\hat{\pi}_{t} 
    \,dt
$$
\begin{equation}
                                              \label{10.29.5}
+\big[ D_{i}(\phi^{n}_{t}a^{ij}_{t}D_{j}\hat{\pi}_{t})-
(\hat{a}^{ij}_{t}+a^{ij}_{t})(D_{i}\phi^{n}_{t})D_{j}\hat{\pi}_{t}
+ D_{i}(\hat{\pi}_{t}\alpha^{ij}_{t}D_{j}\phi^{n}_{t})
 \big]\,dt,\quad t\leq \tau_{m} .
\end{equation}

As above we conclude that $\phi^{n}_{t}
\hat{\pi}\in\cW^{1}_{p}(\tau_{m})$ 
and that, owing to \cite{Kr09_3}, with probability
one $\phi^{n}_{t}\hat{\pi}_{t}$ is a continuous $\cL_{p}$-valued
process and (a.s.)
$$
\|\phi^{n}_{t\wedge\tau_{m}}
\hat{\pi}_{t\wedge\tau_{m}}\|_{\cL_{p}}^{p}
=\|\phi^{n}\hat{\pi}_{0}\|_{\cL_{p}}^{p}
 + I^{1n}_{t}+I^{2n}_{t} +I^{3n}_{t}  
$$
for all $t$, where
$$
I^{1n}_{t}=-p(p-1)\int_{0}^{t\wedge\tau_{m}}
a^{ij}_{s}\int_{\bR^{d}}
|\phi^{n}_{s}|^{p}|\hat{\pi}_{s}|^{p-2}(D_{i}\hat{\pi}_{s})
D_{j}\hat{\pi}_{s}\,dx\,ds\leq0,
$$
$$
I^{2n}_{t}= -\int_{0}^{t\wedge\tau_{m}} 
\big[D_{i}\hat{b}^{i}_{s}- D_{i}b^{i}_{s}\big]
\int_{\bR^{d}}
|\phi^{n}_{s}|^{p}|\hat{\pi}_{s}|^{p}
\,dx\,ds,
$$
$$
I^{3n}_{t}=\int_{0}^{t\wedge\tau_{m}}
\int_{\bR^{d}}
 |\hat{\pi}_{s}|^{p}\psi^{n}_{s}\,dx\,ds,
$$   
$$
\psi_{s}^{n}=
 pa^{ij}_{s}D_{ij}|\phi^{n}_{s}
|^{p}+ (p-1)(p-2)|\phi^{n}_{s}|^{p-2}\alpha^{ij}_{s}
(D_{i}\phi^{n}_{s})D_{j}\phi^{n}_{s}
$$
$$
+(2-p)|\phi^{n}_{s}|^{p-1}
\alpha^{ij}_{s}D_{ij}\phi^{n}_{s}\big]
+(b^{i}_{s}-\hat{b}^{i}_{s})D_{i}|\phi^{n}_{s}|^{p} ,
$$
where for simplicity of notation the argument $x$
is dropped.

Observe that $|D\phi^{n}_{s}|
\leq 1/n$ and
for $s\leq\tau_{m}$ we have
$|b _{s}-\hat{b} _{s}|\leq N(1+|x|)$, where $N$ is 
independent of $s,x$, and $\omega$.
Furthermore, $D\phi^{n}_{s}\to0$ as $n\to\infty$ and
for $s<\tau_{m}$
$$
 |x|\,|D \phi^{n}_{s}(x)|=
\frac{|x|}{n}\big|(D \phi)\big(\frac{x-\xi_{s}}{n}\big)\big|
\leq \frac{|\xi_{s}|}{n}  
+\frac{|x-\xi_{s}|}{n}\,
\big|(D \phi)\big(\frac{x-\xi_{s}}{n}\big)\big|
$$
$$
\leq m+\sup_{y}|y|\,|D\phi(y)|.
$$
  By adding that
$\hat{\pi}\in\bW^{1}_{p}(T)$, we conclude that
$I^{3n}_{t}\to 0$ uniformly in $t$ (a.s.).
Analyzing $I^{1n}_{t}$ and $I^{2n}_{t}$ is almost trivial
and 
$$\|\phi^{n}_{t\wedge\tau_{m}}
\hat{\pi}_{t\wedge\tau_{m}}\|_{\cL_{p}}^{p}
\to \| \hat{\pi}_{t\wedge\tau_{m}}\|_{\cL_{p}}^{p}
$$
as $n\to\infty$ by the monotone convergence theorem.
It follows that (a.s.) for all $t$
$$
\| 
\hat{\pi}_{t\wedge\tau_{m}}\|_{\cL_{p}}^{p}
=\| \hat{\pi}_{0}\|_{\cL_{p}}^{p}
 -\int_{0}^{t\wedge\tau_{m}} 
(D_{i}\hat{b}^{i}_{s}- D_{i}b^{i}_{s})
 \|\hat{\pi}_{s}\|_{\cL_{p}}^{p} \,ds
$$
$$
 -p(p-1)\int_{0}^{t\wedge\tau_{m}}
a^{ij}_{s}\int_{\bR^{d}}
 |\hat{\pi}_{s}|^{p-2}(D_{i}\hat{\pi}_{s})
D_{j}\hat{\pi}_{s}\,dx\,ds .
$$

Obviously on can drop $\tau_{m}$ in this formula
and then obtain that (a.s.) for all $t\leq T$
$$
G_{t} \| 
\hat{\pi}_{t }\|_{\cL_{p}}^{p}
=\| \hat{\pi}_{0}\|_{\cL_{p}}^{p}
 -p(p-1)\int_{0}^{t }G _{s}
a^{ij}_{s}\int_{\bR^{d}}
 |\hat{\pi}_{s}|^{p-2}(D_{i}\hat{\pi}_{s})
D_{j}\hat{\pi}_{s}\,dx\,ds,
$$
which implies that $G_{t} \| 
\hat{\pi}_{t }\|_{\cL_{p}}^{p}$ is decreasing and continuous
(a.s.).  
 Furthermore, 
since $\phi^{n}_{t}\hat{\pi}_{t}$ are continuous $\cL_{p}$-valued
processes, $\hat{\pi}_{t}$ is at least a weakly continuous
$\cL_{p}$-valued function, but since
$\|\hat{\pi}_{t}\|^{p}_{\cL_{p}}$ is (absolutely)
continuous, $\hat{\pi}_{t}$ is   strongly continuous.
   This proves the lemma.
 \begin{remark}
After we know that 
$\hat{\pi}_{t}$ is a continuous
$\cL_{p}$-valued process on $[0,T]$
the last assertion of Lemma \ref{lemma 10.30.1}
can be also obtained from uniqueness of solutions
of \eqref{9.29.01}
because the $\tau $ in Lemma \ref{lemma 10.30.1}
 is a stopping time and $\hat{\pi}_{t\wedge\tau}$
is obviously a solution of \eqref{9.29.01} implying that
on the set where $\tau<T$ we have $\hat{\pi}_{t}=0$
for $\tau\leq t\leq T$.
\end{remark}

Before stating the following lemma we introduce a stipulation
accepted throughout the rest of the paper
that if we are given a function $\xi(t,x,y)$, then
we denote 
\begin{equation}
                                                    \label{10.10.1}
 \tilde{\xi}_{t}=\xi_{t}(x_{t})
=\xi(t,x_{t},y_{t}).
\end{equation}
The reader encountered
above already one of these abbreviated notation
(see \eqref{1.16.1}).  
\begin{lemma}
                                        \label{lemma 10.29.1}
Introduce 
$$
\tilde{w} _{t}=\int_{0}^{t}\Psi _{s}\Theta_{s}\,dw _{s},
 \quad 
\tilde{\sfB}_{t}=\sfB_{t}(x_{t})=\Psi(t,y_{t})B(t,x_{t},y_{t}).
$$
Then $\tilde{w} _{t}$ is a Wiener process and
the process
$$
\rho_{t}=\rho_{t}(\tilde{\sfB},d\tilde{w})
 =\exp(-\int_{0}^{t}
\tilde{\sfB}^{k}_{s}\,d\tilde{w}^{k}_{s}
-\tfrac{1}{2}\int_{0}^{t}|\tilde{\sfB}_{s}|_{\ell_{2}}^{2}\,ds)
$$ 
is a martingale on $[0,T]$. 
\end{lemma}

Proof.
The first assertion follows from L\'evy's theorem.
To prove the second one observe that
$$
\int_{0}^{t}
\tilde{\sfB}^{k}_{s}\,d\tilde{w}^{k}_{s}=
\int_{0}^{t}
\tilde{\sfB}^{*}_{s}\Psi _{s}\Theta_{s}\,dw _{s}. 
$$
Furthermore, the system
$$
dx_{t}=\big(b(t,z_{t})-\theta(t,y_{t})
\Theta^{*}(t,y_{t})\Psi^{2}(t,y_{t})
 B(t,z_{t})\big)\,dt+\theta(t,y_{t})\,dw_{t},
$$
$$
dy_{t}=\Theta(t,y_{t})\,dw_{t},
$$
which is obtained from \eqref{eq3.2.14} by formal
application of the measure change, has a unique solution
with initial data $z_{0}$ since its coefficients
are locally Lipschitz in $z$ and grow as $|z|\to\infty$
not faster than linearly.
In this situation by the Liptser-Shiryaev theorem
  $\rho$ is a martingale since
$$
\int_{0}^{T}| 
\Psi (t,y(t))B(t,x(t),y(t))|^{2}\,dt<\infty
$$
for any deterministic  functions
$x(t)$ and $y(t)$ which are continuous on $[0,T]$.
The lemma is proved.

\begin{lemma}
                                        \label{lemma 10.8.1}
The process $\hat{\pi}_{t}$ is $\cF^{y}_{t}$-adapted.
\end{lemma}

Proof. Observe that in equation \eqref{9.29.01} we have
$$
d\tilde{y}^{k}_{t}=\Psi^{kr}_{t}\,dy^{r}_{t}=d\tilde{w}^{k} _{t}
+\tilde{\sfB}^{k}_{t}\,dt,
$$
where, as it is pointed out above,
$\tilde{w} _{t}$ is a Wiener process. Furthermore, the processes
$\tilde{y}^{k}_{t}$
is $\cF^{y}_{t}$-adapted since such are $\Psi^{kr}_{t}$ and
equation \eqref{9.29.01} is rewritten as
$$
d\hat{\pi}_{t}=
-\sigma^{ik}_{t}D_{i}\hat{\pi}_{t}\,d\tilde{w}^{k}_{t}
+\big[ D_{i}( a^{ij}_{t}
D_{ j} \hat{\pi }_{t})-b^{i}_{t}D_{i} \hat{\pi}_{t} 
$$
\begin{equation}
                                                  \label{10.7.1}
+( \hat{b}^{i}_{t}-\sigma^{ik}_{t}\tilde{\sfB}^{k}_{t}
)D_{i}\hat{\pi}_{t}   \big]\,dt
,\quad t\leq T .
\end{equation}
Here 
$\sigma^{ik}_{t}\tilde{\sfB}^{k}_{t}$ is independent of 
$x$ and for each $\omega$ the trajectories of
$\sigma^{ik}_{t}\tilde{\sfB}^{k}_{t}$ are locally 
bounded on $\bR_{+}$, which
shows that in order to be able to apply
Theorem \ref{theorem 10.6.1} it only remains to  
refer to Lemma \ref{lemma 10.29.1}.
 The lemma is proved.

\begin{lemma}
                                       \label{lemma 10.8.2}
The
assertions (i)-(iii) of Theorem \ref{theorem 11.8.1}
hold for $\bar{\pi}_{t}:=e^{-Q_{t}}\hat{\pi}_{t}$.
 
\end{lemma}

Proof. Assertion (i) of Theorem \ref{theorem 11.8.1}
follows immediately from Lemma 
\ref{lemma 10.30.1}, the continuity
of $Q_{t}$, and the boundedness of $W_{t}=(D_{ij}Q_{t})$
 away from zero.

To prove  assertion (ii) notice
that $\hat{\pi}\in\bW^{1}_{p}(T)$ and
$$
\int_{0}^{t}\|\bar{\pi}_{s}\|_{W^{1}_{p}}^{p}\,ds
$$
is an $\cF^{y}_{t}$-adapted continuous process on
$[0,T]$. Then after
introducing
$$
\tau'_{m}=T\wedge\inf\{t\geq0: 
\int_{0}^{t}\|\bar{\pi}_{s}\|_{W^{1}_{p}}^{p}\,ds\geq m\}
$$
we get that $\bar{\pi}\in\bW^{1}_{p}(\tau'_{m})$
and $\tau'_{m}=T$ for all large $m$ (a.s.).

We now prove that $\bar{\pi}$ satisfied \eqref{10.20.2}  
define 
$\Phi_{t}=\Psi_{t}^{-1}$
and observe that 
$$
(d\tilde{y}^{k}_{t})d\tilde{y}^{r}_{t}=\delta^{kr}\,dt,\quad
dy^{k}_{t}=\Phi^{kr}_{t}\,d\tilde{w}^{r}_{t}+\tilde{B}^{k}_{t}\,dt,
$$
( $\tilde{B}_{t}=B(t,z_{t})$).
Recall that $\eta_{t}(x)=\exp(-Q_{t}(x))$
satisfies
equation \eqref{9.10.5} for each $x$ with
probability one for all $t\in[0,T]$. It turns out that
this equation also holds in the sense of generalized
functions. Owing to the special structure
of $Q_{t}$, this follows from the stochastic version
of Fubini's theorem (see, for instance,
Lemma 2.7  of \cite{Kr09_4}).  

Next, for $m=1,2,...$ set
\begin{equation}
                                               \label{11.2.1}
\tau''_{m}=T\wedge\inf\{t\geq0:|z_{t}|+|DQ_{t}(0)|\geq m\}.
\end{equation}
Note that for a constant $N_{0}$ independent
of $m$ for $t<\tau''_{m}$ we have
$$
|\sfB_{t}(x)| +|b_{t}(x)|\leq N_{0}(1+|x|+m),
\quad|\tilde{\sfB}_{t}|+|\tilde{B}_{t}|
\leq N_{0}(1+2m).  
$$
Furthermore, $D_{i}Q_{t}(x)=x^{j}D_{ij}Q_{t}+D_{i}Q_{t}(0)$,
so that increasing $N_{0}$ if needed we may assume that
for $t<\tau''_{m}$
$$
|DQ_{t}(x)|\leq N_{0}(1+ |x|+m).
$$
Then as is easy to see (cf. \eqref{11.1.1})
$u_{t}:=\hat{\pi}_{t}$ and 
$\tilde{u}_{t}:=\eta_{t}$
satisfy the condition of
Theorem~\ref{theorem 10.14.1} with appropriate
$f,\tilde{f},g,\tilde{g}$ and $\tau''_{m}$
in place of $\tau$.

By Theorem \ref{theorem 10.14.1}
in the sense of generalized functions
$$
d(\eta_{t}\hat{\pi}_{t})
=I^{r}_{t}\,
d\tilde{y}^{r}_{t}+J_{t}\,dt,\quad t\leq\tau''_{m},
$$
where
$$
I^{r}_{t}=\hat{\pi}_{t}\Lambda^{r*}_{t}\eta_{t}  
-\eta_{t} \sigma^{ir}_{t}D_{i}\hat{\pi}_{t}
=\Lambda^{r*}_{t}(\eta_{t}
\hat{\pi}_{t}),
$$
$$
J_{t}=-(\eta_{t} \sfB^{k}_{t}-
\sigma^{ik}_{t}D_{i}\eta_{t} 
)\sigma^{jk}_{t}D_{j}\hat{\pi}_{t} +\hat{\pi}_{t}
  L^{*}_{t}  \eta_{t}
$$
$$ 
+\eta_{t}\big[  a^{ij} _{t}D_{ij}
 \hat{\pi }_{t}-b^{i}_{t}D_{i}\hat{\pi}_{t}
 +(  \sigma^{ik}_{t}\sfB^{k}_{t}+
2\eta^{-1}\hat{a}^{ij}_{t} 
D_{j}\eta_{t})D_{i}\hat{\pi}_{t} 
 \big] 
$$
$$
=\hat{\pi}_{t}
  L^{*}_{t}  \eta_{t} +\eta_{t}(a^{ij} _{t}D_{ij}
 \hat{\pi }_{t}-b^{i}_{t}D_{i}\hat{\pi}_{t})
+2a^{ij}_{t}(D_{i}\hat{\pi}_{t})
D_{j}\eta_{t} =L^{*}_{t}(\eta_{t}\hat{\pi}_{t}).
$$
 In other words (see Theorem \ref{theorem 10.14.1})
for any $\phi\in C^{\infty}_{0}$
with probability one
$$
(\bar{\pi}_{t\wedge\tau''_{m}},\phi)
=(\bar{\pi}_{ 0},\phi)
+\int_{0}^{t}I_{s\leq\tau''_{m}}(\bar{\pi}_{s},\Lambda^{k}
_{s}\phi)\,d\tilde{y}^{k}_{s}
+\int_{0}^{t}I_{s\leq\tau''_{m}}(\bar{\pi}_{s}, L_{s}
\phi)\,ds
$$
for all $t\geq0$. Obviously, one can take here $\tau_{m}:
=\tau'_{m}\wedge\tau''_{m}$ in place of $\tau''_{m}$
and then after recalling that $\bar{\pi}\in\bW^{1}_{p}
(\tau'_{m})$ one concludes that $\bar{\pi}$
is a solution of \eqref{10.20.2}
  in the sense of Definition \ref{definition 3.20.01}.  
The final assertion in (iii) is obtained in the same way
as \eqref{11.1.1}. The  lemma is proved.

To better orient the reader it is worth noting that 
in the next lemma the second factor on the left
in \eqref{11.5.1} contains the negative of two terms 
in \eqref{11.1.3}.

\begin{lemma}
                                    \label{lemma 11.5.1}
We have
$$
\rho_{t}(\tilde{\sfB},d\tilde{w})
\exp\big(-\int_{0}^{t}(V^{*}_{s}W^{-1}_{s}\dot{\sfB}_{s}
-\sfB^{*}_{s}(0))\,d\tilde{y}_{s}
-\tfrac{1}{2}\int_{0}^{t}|\dot{\sfB}^{*}_{s}W_{s}^{-1}
V_{s}-\sfB_{s}(0)|^{2}\,ds\big)
$$
\begin{equation}
                                                \label{11.5.1}
=\rho_{t}(\tilde{\sfB}- \sfB_{\cdot}(0)+
\dot{\sfB}^{*}W^{-1}V,d\tilde{w}).
\end{equation}
Furthermore, the right-hand side is a martingale on $[0,T]$.
\end{lemma}

Proof. The equality is obtained by
simple manipulations. As in the proof of Lemma
\ref{lemma 10.29.1}, to prove that \eqref{11.5.1} is
a martingale we are going to use the
Liptser-Shiryaev theorem by considering the system
consisting of \eqref{9.29.4}, \eqref{9.12.1},  and
\eqref{9.12.2}. We do not include
\eqref{9.15.1}  because  
  $U_{t}$ does not enter \eqref{11.5.1}.  First of all we find
a smooth bounded, uniformly
nondegenerate $d
\times d$-matrix-valued function
  $F(W )$  such that   $F(W_{t})=W_{t}$. The fact that
this is possible follows from
Lemma \ref{lemma 11.1.1}. Then set
$$
A(t,z,W,V)=\Theta^{*}(t,y)
 \Psi^{2}(t,y )
\big(B(t,z )-B(t,0,y )+\dot{B}^{*}(t,y)F^{-1}(W )
V \big).
$$
 
After changing the probability measure formally
we arrive at the   system consisting of \eqref{9.12.1}
with $\sigma_{t}=\sigma(t,y_{t})$, $\hat{a}_{t}
=\hat{a}(t,y_{t})$, and with $F(W_{t})$
in place of $W_{t}$ on the right and the following two equations
$$
dz_{t}=\check{\theta}(t,y_{t})\,dw_{t}
+\big[\check{b}(t,z_{t})-\check{\theta}(t,y_{t})
A(t,z_{t},W_{t},V_{t})\big]\,dt ,
$$
$$
dV_{t}=-\big(F(W_{t})\sigma(t,y_{t})+\dot{B}(t,y_{t})
\Psi(t,y_{t})\big)\Psi(t,y_{t})
\Theta(t,y_{t})\,dw_{t}
$$
$$
+\big(F(W_{t})\sigma(t,y_{t})+\dot{B}(t,y_{t})
\Psi(t,y_{t})\big)\Psi(t,y_{t})
\Theta(t,y_{t})A(t,z_{t},W_{t},V_{t})\,dt
$$
$$
-\big(F(W_{t})\sigma(t,y_{t})+\dot{B}(t,y_{t})
\Psi(t,y_{t})\big)\Psi(t,y_{t})B(t,z_{t})\,dt
$$
$$
+\big[\big(\dot{B}(t,y_{t})
\Psi(t,y_{t})\sigma^{*}(t,y_{t})-\dot{b}(t,y_{t})\big)V_{t}
-2F(W_{t})\hat{a}(t,y_{t})V_{t}\big]\,dt
$$
$$
+\big[F(W_{t})\big(\sigma(t,y_{t})\Psi(t,y_{t})B(t,0,y_{t})-
b (t,0,y_{t}) \big)
+\dot{B}(t,y_{t})
\Psi^{2}(t,y_{t})B (t,0,y_{t})
\big]\,dt.
$$
 
This system has a unique solution with prescribed
initial data since its coefficients are locally Lipschitz continuous
and may grow to infinity as $|z|+|W|+|V|\to\infty$ not faster
than linearly.
Moreover,
$$
\int_{0}^{T}|A(t,z(t),W(t),V(t))|^{2}\,dt<\infty
$$
for any functions $z(t),W(t),V(t)$ which are continuous
on $[0,T]$. 
This implies that 
process \eqref{11.5.1} is a martingale on $[0,T]$ and
the lemma is proved.

{\bf Proof of Lemma \ref{lemma 11.1.1}}.
Notice that \eqref{9.15.1} yields $U_{t}$
once   $W_{t}$ and $V_{t}$ are found.
Equation \eqref{9.12.2} is linear with respect to
$V_{t}$ and proving the existence and uniqueness
of its solution presents
no difficulty if $W_{t}$ is known.

Equation \eqref{9.12.1} can be considered for each $\omega$
separately. Then the theory of ODEs allows us to
conclude that a unique solution exists until it 
blows up and it is  $\cF^{y}_{t}$-adapted.
Uniqueness  implies that $W_{t}=W^{*}_{t}$.
Furthermore, at least on a small time interval $W_{t}>0$.
It turns out that $W_{t}>0$ on any interval of time where
$W_{t}$ is bounded.

Indeed, if not, then for some   $t_{0}>0 $ we would have that
$\det W_{t_{0}}=0$,  $W_{t}$ is bounded on $[0,t_{0}]$
and $\det W_{t}>0$ for $t<t_{0}$.
However, for $t<t_{0}$
\begin{equation}
                                               \label{11.8.1}
\frac{d}{dt}\,\det W_{t }=\tr \dot{W}_{t}W^{-1}_{t}\det W_{t },
\end{equation}
and 
$$
\tr \dot{W}_{t}W^{-1}_{t}=2\tr(\dot{\sfB}_{t}\sigma_{t}^{*}
-\dot{b}_{t} )-2\tr \hat{a}_{t}W_{t}+\tr\dot{\sfB}_{t}\dot{\sfB}
^{*}_{t}W^{-1}_{t},
$$
where the last term is nonnegative as the trace of the product
of two symmetric nonnegative matrices. It follows, that
$\tr \dot{W}_{t}W^{-1}_{t}$ is bounded from below on $[0,t_{0})$
and hence equation \eqref{11.8.1} implies that $\det W_{t_{0}}
>0$.

Next, it turns out that the solution does not blow up on $[0,T]$.
Indeed
$$
\frac{d}{dt}\,\tr W_{t}W _{t}=4\tr(\dot{\sfB}_{t}\sigma_{t}^{*}
-\dot{b}_{t} )W_{t}W _{t}+ 2\tr\dot{\sfB}_{t}
\dot{\sfB}_{t}^{*}W_{t}-4\tr \hat{a}_{t}
W^{3}_{t}  ,
$$
where the last trace is  nonnegative again on the interval 
of existence of $W_{t}$. 
Here 
$$
\tr\dot{\sfB}_{t}
\dot{\sfB}_{t}^{*}W_{t}\leq N(\tr W^{2}_{t})^{1/2}
\leq N+\tr W^{2}_{t},
$$
where $N$ is a constant. Also for two    matrices
$A$ and $W$ such that $W$ is symmetric and nonnegative it holds that 
$$
(\tr AW^{2})^{2}\leq\|A\|\,\|W^{2}\|\leq
\|A\|(\tr W^{2})^{2}.
$$
This and Gronwall's inequality
imply that $W_{t}$ is bounded on $[0,T]$.  
Obviously the bound of $W_{t}$ is uniform with respect to
$\omega$. The lower bound is also uniform since by the above
$\det W_{t}$ is bounded away from zero on $[0,T]$ uniformly
with respect to $\omega$. The lemma is proved.

\mysection{Proof of Theorems \protect\ref{theorem 11.8.1}
and \protect\ref{theorem 11.11.1}}

                                           \label{section 1.1.4}
Take a function $\varphi\in C^{\infty}_{0}(\bR^{d_{1}})$
and let $c (t,y)$
be a smooth, bounded, and nonnegative function on 
$[0,T]\times\bR^{d_{1}-d}$. Recall that 
the operator $\check{L}$ is introduced in \eqref{eq3.2.19.1}
and consider the following deterministic problem
$$
\partial_{t}v (t,z)+
\check{L} v (t,z)-c (t,y)v (t,z) =0,\quad
t\in[0,T],z\in\bR^{d_{1}}, 
$$
\begin{equation} 
                                          \label{10.13.1}
v (T,z) =\varphi(z),\quad z\in\bR^{d_{1}}.
 \end{equation}
\begin{remark}
                                     \label{remark 9.20.1}

By Theorem 2.5 of \cite{KP}, for any $\alpha\in(0,1)$
there exists a unique classical solution
$v $ of \eqref{10.13.1} such that, for any $t\in[0,T]$,
$v(t,\cdot)\in C^{2+\alpha}(\bR^{d_{1}})$ and the standard
$C^{2+\alpha}(\bR^{d_{1}})$-norms of $v(t,\cdot)$ are bounded
on $[0,T]$. 
If we denote by $z_{t}(s,z)$, $t\geq s$, the solution of 
system \eqref{eq3.2.14} which starts at $z$ at moment $s\leq T$,
then by It\^o's formula we have  
$$
v (s, z)=E\varphi(z_{T }(s,z))\exp(-\int_{s}^{T }
c_{r}( y_{r}(s,z))\,dr),
$$
$$
|v(s,z)|\leq \sup|\varphi|P\{\tau(s,z)\leq   T\}
\leq \sup|\varphi|e^{N_{0}T}E e^{-N_{0}\tau(s,z)},
$$
where $N_{0}>0$
is an arbitrary constant,
$\tau(s,z)$ is the first time $z_{t}(s,z)$ hits
$\{z:|z|\leq R\}$, and $R$ is such that $\varphi(z)=0$
for $|z|\geq R$.  Take an $m\geq0$ and introduce
$\psi(z)=(1+|z|^{2})^{-m}$. It is not hard to see that,
if $N_{0}$ is sufficiently large, then
$$
\check{L}_{t}\psi(z)-N_{0}\psi(z) \leq0.
$$
By It\^o's formula, for $|z|\geq R$,
$$
\psi(R)E e^{-N_{0}\tau(s,z)}\leq \psi(z),
$$
implying that for any $m\geq0$ there is a constant
$N$ such that for all $(s,z)$
$$
|v(s, z)|\leq\frac{N}{(1+|z|^{2})^{m}}.
$$

The argument in the proof of Lemma 4.11 of \cite{KZ}
proves that the same estimate holds for $\partial v(s,z)/
\partial z^{i}$
and $\partial^{2} v(s,z)/
\partial z^{i}\partial z^{j}$, $i,j=1,...,d_{1}$.
\end{remark}

Before we come to a crucial point we state the following.
\begin{lemma}
                                      \label{lemma 10.29.2}
Let $\xi_{t}$ be a nonnegative
continuous  martingale on $[0,T]$
and let $\zeta_{t}$ be a continuous $\cF_{t}$-adapted process
given on $[0,T]$
such that $\xi_{t}\zeta_{t}$ is a local martingale on
$[0,T)$. Assume that
$$
E\xi_{T}\sup_{[0,T]}|\zeta_{t}|<\infty.
$$
Then $\xi_{t}\zeta_{t}$ is a   martingale on
$[0,T]$.
\end{lemma}

Proof. We need to prove that for any stopping time
$\tau\leq T$ we have $E\xi_{\tau}\zeta_{\tau}=E\xi_{0}\zeta_{0}$.
Here the left hand side equals $E\xi_{T}\zeta_{\tau}$
and we are given that there exists a sequence of stopping times 
$\tau_{n}\uparrow T$ such that $E\xi_{T}\zeta_{\tau
\wedge\tau_{n}}=E\xi_{0}\zeta_{0}$. Using the dominated 
convergence theorem yields the desired result 
and proves the lemma.

\begin{lemma}
                                       \label{lemma 10.8.02}
The process
$$
  \rho_{t}e^{-\int_{0}^{t}c_{s}(y_{s})\,ds}\int_{\bR^{d}}
v (t,x,y_{t})\bar{\pi}_{t}(x)\,dx
$$
is a martingale on $[0,T]$.
\end {lemma}

Proof. Define ($c_{t}=c (t,y_{t})$, $v_{t}(x)=v(t,x,y_{t})$)
$$
D^{y}_{k}=\frac{\partial}{\partial y^{k}},\quad
D^{y}_{kr}=D^{y}_{k}D^{y}_{r},\quad
C_{t}=\exp(
-\int_{0}^{t}c_{s}\,ds),\quad
\chi_{t}=C_{t}v_{t}\bar{\pi}_{t}.
$$
We need to show that
\begin{equation}
                                                 \label{11.11.6}
  \rho_{t} \int_{\bR^{d}}
\chi_{t}(x)\,dx
\end{equation}
is a martingale.

Observe that by It\^o's formula 
 and \eqref{10.13.1} we have
$$
d [v_{t}(x )C_{t}]=d [v_{t} C_{t}]=C_{t}\big[
D^{y}_{k}v_{t} \,dy^{k}_{t}
+
\big(\partial_{t}v_{t} -c_{t}v_{t} 
+\check{a}^{kr}_{t}
 D^{y}_{kr} v_{t} \big)\,dt\big]
$$
\begin{equation}
                                   \label{10.24.2}
=C_{t}\big[D^{y}_{k}v_{t} 
\Phi^{kr}_{t}\,d\tilde{w}^{r}_{t}-\big(L_{t}
   v_{t} +2 \check{a}^{ik}_{t}
 D _{i}D^{y}_{k} v_{t}+
 (B^{k}_{t}-\tilde{B}^{k}_{t})
D^{y}_{k}v_{t} \big)\,dt
 \big],
\end{equation}
where    we dropped the arguments $x$
  for shortness and, of course, $D^{y}_{k}v_{t}=
(D^{y}_{k}v) (t,x,y_{t})$,
$D^{y}_{kr}v_{t}=
(D^{y}_{kr}v) (t,x,y_{t})$, and $D_{i}D^{y}_{k}v_{t}=
(D_{i}D^{y}_{k}v) (t,x,y_{t})$. By the way, observe that
$$
 \sigma^{ir}_{t}\Phi^{kr}_{t}
=2\check{a}^{ik}_{t},\quad B^{k}_{t}=\Phi^{kr}_{t}
\sfB^{r}_{t}.
$$

Similarly to the proof   
of Lemma \ref{lemma 10.8.2}   
we conclude that \eqref{10.24.2} holds in the sense
of distributions and that Theorem \ref{theorem 10.14.1}
is applicable to
$v_{t}\bar{\pi}_{t}$ on the time interval $t\leq\tau_{m}$
for any $n$, where $\tau_{m}$ are taken from 
Lemma \ref{lemma 10.8.2}.
It follows that for any $m$ for $t\leq\tau_{m}$
$$
d\chi_{t}= C_{t}(\bar{\pi}_{t}\Phi^{kr}_{t}D^{y}_{k}v_{t} 
+v_{t}\Lambda^{r*}_{t}\bar{\pi}_{t})
\,d\tilde{w}^{r}_{t}
$$
$$
-C_{t}\bar{\pi}_{t}\big(L_{t}
   v_{t} + \sigma^{ir}_{t}\Phi^{kr}_{t}
 D _{i}D^{y}_{k} v_{t}+\Phi^{kr}_{t}
 (\sfB^{r}_{t}-\tilde{\sfB}^{r}_{t})
D^{y}_{k}v_{t} \big)\,dt
$$
$$
+C_{t}v_{t}(L^{*}_{t}\bar{\pi}_{t}+
\tilde{\sfB}_{t}^{k}
\Lambda^{k*}_{t}\bar{\pi}_{t})\,dt
+C_{t}(D^{y}_{k}v_{t}) 
\Phi^{kr}_{t}\Lambda^{r*}_{t}\bar{\pi}_{t}\,dt.
$$
It is convenient to rearrange the above terms 
by using the notation
$$
\zeta^{r}_{t}=C_{t}(\bar{\pi}_{t}\Phi^{kr}_{t}D^{y}_{k}v_{t} 
+v_{t}\Lambda^{r*}_{t}\bar{\pi}_{t}).
$$
We have
$$
d\chi_{t}=\zeta_{t}^{r}\,d\tilde{w}^{r}_{t}
+(\tilde{\sfB}^{r}_{t}\zeta^{r}_{t}+I^{1}_{t}+I^{2}_{t})\,dt,
\quad t\leq\tau_{m},
$$
where
$$
I^{1}_{t}= C_{t} (v_{t}L^{*}_{t}\bar{\pi}_{t}
-\bar{\pi}_{t}L_{t} v_{t}),\quad
I^{2}_{t}=-C_{t}\Phi_{t}^{kr}\sigma^{ir}_{t}(
\bar{\pi}_{t}
\sigma^{ir}_{t}D_{i}D^{y}_{k}v_{t}+
(D^{y}_{k}v_{t})D_{i}\bar{\pi}_{t} )
$$
$$
=-C_{t}\Phi_{t}^{kr}\sigma^{ir}_{t}D_{i}(\bar{\pi}_{t}
D^{y}_{k}v_{t}).
$$
In the integral form this means that for any $\phi
\in C^{\infty}_{0}$ with probability one 
$$
(\chi_{t\wedge\tau_{m}},\phi)=(\chi_{0},\phi)
+\int_{0}^{t}I_{s\leq\tau_{m}} (\zeta^{r}_{s},\phi)\,
d\tilde{w}_{s}^{r}+
\int_{0}^{t}I_{s\leq\tau_{m}}
\tilde{\sfB}^{r}_{s}(\zeta^{r}_{s},\phi)\,ds
$$
$$
+\int_{0}^{t}I_{s\leq\tau_{m}}C_{s}a^{ij}_{s}
(\bar{\pi}_{s}D_{j}v_{s}-v_{s}D_{j}\bar{\pi}_{s},
D_{i}\phi)\,ds
$$
\begin{equation}
                                     \label{10.24.4}
+\int_{0}^{t}I_{s\leq\tau_{m}}C_{s}
\big[(\bar{\pi}_{s}v_{s},b^{i}_{s}D_{i} \phi )+
\Phi^{kr}_{s}\sigma^{ir}_{s}(\bar{\pi}_{t}D^{y}_{k}v_{s},
D_{i}\phi)\big]\,ds.
\end{equation}
 We take a $\phi$ such that $\phi(0)=1$ and plug 
$\phi_{j}$ into
\eqref{10.24.4} in place of $\phi$, where
$\phi_{j}(x)=\phi(x/j)$, $j=1,2,...$.

Observe that 
$$
(\zeta^{r}_{s},\phi_{j})=C_{s}(\phi_{j} ,\bar{\pi}_{s}
 \Phi^{kr}_{s}D^{y}_{k}v_{s}+v_{s}
\Lambda^{r*}_{s}\bar{\pi}_{s})
$$
and for any $r$ and $k$
$$
\int_{0}^{T}(1,|\bar{\pi}_{s}D_{k}^{y}v_{s}|+
|v_{s}\Lambda^{r*}_{s}\bar{\pi}_{s}|)^{2}\,ds
$$
$$
\leq N\int_{0}^{T}\|\bar{\pi}_{s}\|_{W^{1}_{p}}^{2}
\|v_{s}\|_{W^{1}_{p'}}^{2}\,ds
\leq N\|\bar{\pi}\|_{\bW^{1}_{p}(T)}^{2}<\infty,
$$
where $N$ is independent of $\omega$.
 By the dominated convergence theorem and the rules
for passing to the limit under the sign of stochastic
integral it follows that in probability
uniformly on $[0,T]$ 
$$
\int_{0}^{t}I_{s\leq\tau_{m}} (\zeta^{r}_{s},\phi_{j})\,
d\tilde{w}_{s}^{r}\to
\int_{0}^{t}I_{s\leq\tau_{m}} 
C_{s}(1,\bar{\pi}_{s}
 \Phi^{kr}_{s}D^{y}_{k}v_{s}+v_{s}
\Lambda^{r*}_{s}\bar{\pi}_{s})\,d\tilde{w}_{s}.
$$
Similarly, and in an easier fashion one analyzes the remaining
terms in \eqref{10.24.4} and concludes that
for any $m$  
$$
d(\chi_{t},1)=
C_{t}(1,\bar{\pi}_{t}
 \Phi^{kr}_{t}D^{y}_{k}v_{t}+v_{t}
\Lambda^{r*}_{t}\bar{\pi}_{t})\,d\tilde{y}_{t},\quad t\leq \tau_{m}.
$$
By using It\^o's formula we then
immediately obtain that the process \eqref{11.11.6}
is at least a local martingale on $[0,T]$. We rewrite it
as $\xi_{t}\zeta_{t}$, where (see
Remark \ref{remark 1.17.1} and  Lemma \ref{lemma 11.5.1})
$ \xi_{t}=
\rho_{t}(\tilde{\sfB}- \sfB_{\cdot}(0)+
\dot{\sfB}^{*}W^{-1}V,d\tilde{w})  
$ 
and 
$$
\zeta_{t}=
e^{-A_{t}-\int_{0}^{t}c_{s}(y_{s})\,ds}
\int_{\bR^{d}}\hat{\pi}_{t}(x)
v_{t}(x)\exp\big(-\tfrac{1}{2}\int_{0}^{t}
|W^{1/2}_{ s}x+W^{-1/2}_{s}V_{s}|^{2}\,ds\big)\,dx.
$$
Owing to Lemma \ref{lemma 10.30.1} the process $\zeta_{t}$
is bounded on $[0,T]$ by a constant times
$\|\pi_{0}\|_{\cL_{p}}$ which along with Lemma
\ref{lemma 10.29.2} implies that $\xi_{t}\zeta_{t}$
is a martingale. The lemma is proved.

{\bf Proof of Theorem \ref{theorem 11.8.1}}.
Recall that assertions (i)-(iii) are proved
in Lemma \ref{lemma 10.8.2}.
By Lemma \ref{lemma 10.8.02} and  It\^o's formula
$$
Ee^{-\int_{0}^{T}c_{s}(y_{s})\,ds}\varphi(z_{T})=
Ev(0,x_{0},y_{0})=E\int_{\bR^{d}}v(0,x,y_{0})
\bar{\pi}_{0}\,dx
$$
$$
=E\rho_{T}e^{-\int_{0}^{T}c_{s}(y_{s})\,ds}
\int_{\bR^{d}}\varphi(x,y_{T})\bar{\pi}_{T}(x)\,dx
$$
$$
=E\bar{\rho_{T}}e^{-\int_{0}^{T}c_{s}(y_{s})\,ds}
\int_{\bR^{d}}\varphi(x,y_{T})\bar{\pi}_{T}(x)\,dx,
$$
 where $\bar{\rho_{T}}=E(\rho_{T}\mid\cF^{y}_{T})$.
Since the equality between the extreme terms holds
for sufficiently wide class of functions $c$, we get that
$$
E\big(\varphi(z_{T})\mid\cF^{y}_{T}\big)
=\bar{\rho_{T}}
\int_{\bR^{d}}\varphi(x,y_{T})\bar{\pi}_{T}(x)\,dx\quad
\text{(a.s.)}.
$$
The arbitrariness of $\phi$ implies that 
 $\bar{\pi}_{T}\geq0$ (a.s.) and
$$
1=\bar{\rho_{T}}
\int_{\bR^{d}} \bar{\pi}_{T}(x)\,dx,\quad
(1,\bar{\pi}_{T})=\int_{\bR^{d}} \bar{\pi}_{T}(x)\,dx>0,\quad
\bar{\rho_{T}}=(1,\bar{\pi}_{T})^{-1}
$$
(a.s.). It follows that for any Borel $f\geq0$
equation \eqref{10.13.3} holds with $t=T$.

The above argument can be repeated for any $t\leq T$
by taking $t$ in place of $T$. Then we obtain 
\eqref{10.13.3} for any $t$. Furthermore, for any $t$ we will have that
that $\bar{\pi}_{t}\geq0$ and $(1,\bar{\pi}_{t})>0$
(a.s.). Actually, the last two properties hold
with probability one for all $t$ at once since with
probability one $\bar{\pi}_{t}$ is a continuous 
$\cL_{1}$-function by Lemma \ref{lemma 10.8.2}
and by Lemma \ref{lemma 10.30.1}, on the set where $\tau
=\inf\{t\geq0:(1,\bar{\pi}_{t})=0\}<T$, 
we have $ \bar{\pi}_{T}=0$, which 
only happens with probability zero.
The theorem is proved.

{\bf Proof of Theorem \ref{theorem 11.11.1}}. We use part
of notation from the proof of Lemma \ref{lemma 10.8.02}
but this time 
take $\bar{\pi}_{t}=\eta_{t} = e^{-Q_{t}}$. Then
by It\^o's formula and \eqref{9.10.5} we obtain that
for each $x$  
$$
d\chi_{t}=\zeta^{r}_{t}\,(d\tilde{w}^{r}_{t}+
\tilde{\sfB}^{r} \,dt)+
C_{t}(v_{t}L^{*}_{t}\bar{\pi}_{t}-\bar{\pi}_{t}L_{t}v_{t})
\,dt -C_{t}\Phi^{kr}_{t}\sigma^{ir}_{t}
D_{i}(\bar{\pi}_{t}D^{y}_{k}v_{t})\,dt.
$$
By using the stochastic Fubini theorem and 
integrating by parts we see that
$$
d(\chi_{t},1)=(\zeta^{r}_{t},1)\,(d\tilde{w}^{r}_{t}+
\tilde{\sfB}^{r} \,dt)
$$ 
which implies that process \eqref{11.11.6}
is a local martingale on $[0,T]$.
  We rewrite it
as $\xi_{t}\zeta_{t}$, where (see Remark \ref{remark 1.17.1} and
Lemma \ref{lemma 11.5.1})
$ \xi_{t}=
\rho_{t}(\tilde{\sfB}- \sfB_{\cdot}(0)+
\dot{\sfB}^{*}W^{-1}V,d\tilde{w}) 
$ 
and 
$$
\zeta_{t}=
e^{-A_{t}-\int_{0}^{t}c_{s}(y_{s})\,ds}
\int_{\bR^{d}} 
v_{t}(x)\exp\big(-\tfrac{1}{2}\int_{0}^{t}
|W^{1/2}_{ s}x+W^{-1/2}_{s}V_{s}|^{2}\,ds\big)\,dx.
$$
Notice that $\xi_{t}$ is a martingale and $\zeta_{t}$
is obviously bounded. By Lemma \ref{lemma 10.29.2}
we conclude that process \eqref{11.11.6} is a martingale.
 
After that it suffices to repeat the proof of 
Theorem \ref{theorem 11.8.1} dropping unnecessary here
details concerning the fact that $(1,\bar{\pi}_{t})>0$.
The theorem is proved.

\end{document}